\documentclass[11pt]{article}
\usepackage{amssymb}
\usepackage{amsthm}
\usepackage{graphicx}
\usepackage{graphics}
\usepackage{amsthm}

\input pictex.tex
\newcommand{\esp}{\hspace{0.05cm}}

\theoremstyle{definition}

\newtheorem{thm}{Theorem}[section]

\newtheorem{lem}[thm]{Lemma}

\input epsf

\begin{document}

\date{}

\author{Andr\'es Navas}

\title{On centralizers of interval diffeomorphisms
in critical (intermediate) regularity}
\maketitle

\vspace{-0.5cm}

\hfill{\em To the memory of Sergio Plaza Salinas}

\vspace{0.5cm}

\noindent{\bf Abstract.} We extend to the critical (intermediate) regularity several
results concerning rigidity for centralizers and group actions on the interval.

\vspace{0.15cm}


\vspace{0.5cm}

\noindent{\Large {\bf Introduction}}

\vspace{0.3cm}

Group actions on 1-dimensional manifolds is a well-developed subject that takes its source
from the theory of codimension-1 foliations (see \cite{book} for a general panorama). When
these actions are by smooth-enough (namely, $C^2$) diffeomorphisms, the general picture is
well-understood essentially by the classical works of Denjoy, Sacksteder, and Kopell,
among others. The interest in considering actions of lower regularity comes from
different sources (see, for example, \cite{FF,grigorchuk,tsuboi}). It appears that many
interesting phenomena from both the group theoretical and the dynamical viewpoints arise in
{\em intermediate regularity}, that is, for actions by diffeomorphisms of differentiabiity
classes between $C^1$ and $C^2$. This is the main subject of \cite{CJN,DKN,KN,growth},
where several relevant problems have already been settled. Nevertheless, due to technical
reasons, in many cases it was necessary to avoid certain {\em critical regularities} for
which the so far existing arguments do not apply. Despite of this, it is conjectured
that the corresponding rigidity phenomena should still hold in these critical cases. In
this work, we confirm this intuition for centralizers and group actions on the interval
by providing concrete proofs. According to the (methods of) construction of \cite{tsuboi}
for Theorem A below, \cite{KN} for Theorem B, and \cite{CJN} for Theorem C, our results
are optimal (in the H\"older scale). Unfortunately, one of our main arguments does not
apply in the most important context, namely, that of the generalized Denjoy Theorem
in critical regularity, although it provides important evidence for its validity.

\vspace{0.3cm}

\noindent{\bf\em The generalized Kopell Lemma in critical regularity.} Our first result 
is the extension to the critical regularity of \cite[Theorem B]{DKN}. Actually, this 
may be considered as our main result, as all the next ones are based on similar 
ideas but using more involved combinatorial constructions.

\vspace{0.4cm}

\noindent{\bf Theorem A.} {\em Let $\{ I_{i_1,\ldots,i_{d+1}} \!:
(i_1,\ldots,i_{d+1}) \!\in\! \mathbb{Z}^{d+1} \}$ be a family of subintervals of $[0,1]$
that are disposed respecting the lexicographic order. Assume that $f_1,\ldots,f_{d+1}$
are diffeomorphisms such that \esp $f_j (I_{i_1,\ldots,i_{j-1},i_j,i_{j+1},\ldots,i_{d+1}}) =
I_{i_1,\ldots,i_{j-1},1+i_j,i_{j+1},\ldots,i_{d+1}}$, \esp for all \esp $1 \leq j \leq d+1$.
Then $f_1,\ldots,f_d$ cannot be all of class $C^{1+1/d}$ provided that $f_{d+1}$ is
of class $C^{1+\alpha}$ for some $\alpha > 0$ and commutes with $f_1,\ldots,f_d$.}

\vspace{0.4cm}

Compared to \cite[Theorem B]{DKN}, this result is an improvement in what concerns the
hypothesis of regularity for $f_1,\ldots,f_d$. Nevertheless, we impose an extra regularity
assumption for $f_{d+1}$ (in \cite{DKN}, we just require $f_{d+1}$ to be $C^1$). Moreover,
\cite[Theorem B]{DKN} holds in the case of noncommuting maps (see \cite[Exercise 4.1.36]{book}),
and here we strongly use the fact that $f_{d+1}$ commutes with the other $f_i$'s \esp
(although these $f_i$'s are not assumed to commute between them).

\vspace{0.4cm}

\noindent{\bf\em The generalized Kopell Lemma in critical but different regularities.}
Our second result is an extension of Theorem A inspired by \cite{KN}.

\vspace{0.4cm}

\noindent{\bf Theorem B.} {\em Assume we are in the context of Theorem A but each
$f_i$ is of class $C^{1+\alpha_i}$, with $0 < \alpha_i \leq 1$, for non-necessarily
equal values of $\alpha_i$. Then \esp $\alpha_1 + \cdots + \alpha_d < 1$.}

\vspace{0.4cm}

\noindent{\bf\em No smoothing of the Farb-Franks action in the critical regularity.}
Finally, we extend \cite[Theorem A]{CJN} to the critical regularity. The details
for the statement below are provided in \S \ref{FF}.

\vspace{0.42cm}

\noindent{\bf Theorem C.} {\em Farb-Franks' action of $N_d$ is not
topologically semiconjugated to an action by $C^{1+\alpha}$
diffeomorphisms for \esp $\alpha = \frac{2}{d(d-1)}$.}

\vspace{0.42cm}

As in the case of Theorem A, Theorem C above should have a version for different
regularities for the elements in a canonical generating set, the proof of which
should be a combination of the techniques of \S \ref{kopell-2} and \S \ref{FF}
below. Moreover, a combination of the ideas of \cite[\S 3]{CJN} and \cite[\S 3]{KN}
should allow showing that such an extended result is in fact optimal. The
reader will certainly agree in that including all of this would had
artificially overloaded this already very technical article.

\vspace{0.35cm}

\noindent{\em About the proofs.} Roughly, the proof of all the results above proceeds as follows.
Assume that $g$ is a diffeomorphism of the interval that commutes with many diffeomorphisms
$f_1,\ldots,f_k$, and let $h_n = f_{i_n} \circ \cdots \circ f_{i_1}$ be a ``random composition''
of $n$ factors among these $f_i$'s. Taking derivatives in the equality \esp
$g^k = h_n^{-1} \circ g^k \circ h_n,$ \esp we obtain
$$Dg^k(x) = \frac{Dh_n (x)}{Dh_n (h_n^{-1} g^k \esp h_n (x))} \cdot D g^k (h_n(x))
= \frac{Dh_n (x)}{Dh_n (g^k (x))} \cdot D g^k (h_n(x)).$$
As shown in \cite{CJN,DKN,KN}, whenever the regularity is strictly larger than the corresponding
critical one, it is possible to estimate (uniformly on $n$) the value of the {\em distortion} of
$h_n$, that is, an expression of type \esp $Dh_n(x) / Dh_n(y)$ \esp as above. This allows
showing that the derivatives of the iterates of $g$ are uniformly bounded, which is
impossible unless $g$ is trivial. However, as it was already noticed in the
aforementioned  works, this is no longer possible for the critical regularity because of the
failure of convergence of a certain series. The main new idea consists in noticing that despite of
the absence of uniform control for the distortion, elementary estimates show that its growth (in $n$)
is slow (actually, sublinear). By choosing $n = n(k)$ appropriately, this allows showing that
the growth of the derivatives of $g$ is sublinear, which is impossible. This last issue
was cleverly noticed in \cite{PS}, and we reproduce it below (as a lemma) with its
original proof for the convenience of the reader.

\vspace{0.35cm}

\noindent{\bf Lemma [Polterovich-Sodin]}
{\em If $g \!: I \to I$ is a nontrivial $C^1$ diffeomorphism of a bounded interval,
then there exists an infinite increasing sequence of positive integers $k_j$ such that}
$$\max_{x \in I} D g^{k_j} (x) > k_j.$$

\noindent{\bf Proof.} Let $x_0$ be a point such that $g(x_0) \neq x_0$.
Denoting by $J$ the open interval of endpoints $x_0$ and $g(x_0)$, we have
that $J,g^{-1}(J),g^{-2}(J),\ldots$ are two-by-two disjoint. Therefore,
\begin{equation}\label{finita}
\sum_{k \geq 1} \big| g^{-k}(J) \big| < \infty.
\end{equation}
If the conclusion of the lemma does not hold, then there exists $k_0$ such that
for all $k \geq k_0$ and all $y \in I$, we have \esp $Dg^{-k}(y) \geq 1/k$. \esp
Since \esp $\big| g^{-k}(J) \big| = Dg^{-k}(y_k) |J|$ \esp for a certain
$y_k \in J$, this would imply that
$$\sum_{k \geq k_0} \big| g^{-k}(J) \big| \geq \sum_{k\geq k_0} \frac{|J|}{k} = \infty,$$
which contradicts (\ref{finita}). $\hfill\square$

\vspace{0.3cm}

In order to apply our main argument, we will crucially use the following result of Borichev
\cite{borichev}, which extends prior results of Polterovich and Sodin (valid in the
$C^2$ context) to the $C^{1+\alpha}$ category (see also \cite{other}).

\vspace{0.54cm}

\noindent{\bf Theorem [Borichev]} {\em Let $g$ be a $C^{1+\alpha}$ diffeomorphism of
a closed interval $I$ without hyperbolic fixed points, with $0 \!<\! \alpha \!<\! 1$.
Then letting $C_g$ be the $\alpha$-H\"older constant of \esp $\log(Dg)$,
for every $k \geq 0$, we have}
\begin{equation}\label{general}
\max_{x \in I} Dg^k (x) \leq \exp \big( 3 \esp C_{\!g} \esp
|I|^{\alpha} \esp k^{1-\alpha} \big).
\end{equation}

\vspace{0.2cm}

It is important to point out that although this result is not stated in this way
in \cite{borichev}, it readily follows from it (and its proof). Indeed, Theorem 4
of \cite{borichev} only claims that for $I = [0,1]$, one has the estimate
$$\max_{x \in [0,1]} Dg^k (x) \leq \exp \big( A \esp k^{1-\alpha} \big)$$
for a certain constant $A = A(C_{\!g})$. However, a careful reading of the proof therein
shows that $A = 3 \esp C_{\!g}$ actually works. As the reader will notice, having such a
nice quantitative version will be important in our proof. Moreover, the introduction of the
factor $|I|^{\alpha}$ will be also important. This factor comes from an easy renormalization
argument. Indeed, if we denote by $\bar{g}$ the renormalization of $g$ to the unit interval
(more precisely, \esp $\bar{g} := \varphi_I \circ g \circ \varphi_I^{-1}$, \esp where
$\varphi_I$ is the unique orientation-preserving affine homeomorphism sending $I$ into
$[0,1]$), \esp then (the straightened version of) Borichev's theorem (for the unit
interval) yields
\begin{equation}\label{normalization}
\max_{x \in I} Dg^k (x) =
\max_{y \in [0,1]} D \bar{g}^k(y) \leq \exp(3 \esp C_{\bar{g}} \esp k^{1-\alpha}).
\end{equation}
Since \esp $C_{\bar{g}} = C_{g} |I|^{\alpha}$, \esp this implies (\ref{general}).


\section{Proof of the critical generalized Kopell lemma via a random walk argument}
\label{uno}

\hspace{0.45cm} To prove Theorem A, we denote $g:= f_{d+1}$ and we consider a
composition $h_n = f_{i_n} \circ \cdots \circ f_{i_2} \circ f_{i_1}$
of the $f_i$'s. Then for each $k \geq 1$, we have
$$g^k \esp = \esp h_n^{-1} \circ g^k \circ h_n,$$
which yields
$$Dg^k (x)
\esp = \esp \frac{Dh_n (x)}{Dh_n (h_n^{-1} g^k \esp h_n (x))} \cdot D g^k (h_n(x))
\esp = \esp \frac{Dh_n (x)}{Dh_n (g^k (x))} \cdot D g^k (h_n(x)).$$
We will restrict this equality to \esp $x$ in the interval $I$ defined as
the convex closure of \esp $\bigcup_{i_{d+1} \in \mathbb{Z}} I_{0,0,\dots,0,i_{d+1}}$.
\esp (Notice that $I$ is invariant under $g$.) Let $C$ be a simultaneous
$\frac{1}{d}$-H\"older constant for $\log(D f_i)$, where $1 \leq i \leq d$.
Denoting $h_j := f_{i_j} \circ \cdots \circ f_{i_1}$ whenever
$0 \leq j \leq n$ and letting $y = y_k := g^{k}(x) \in I$, we have
\begin{eqnarray*}
\left| \log \left( \frac{Dh_n (x)}{Dh_n (g^k (x))} \right) \right|
&=&
\left| \log \left( \frac{\prod_{j=1}^{n} Df_{i_j} (h_{j-1}(x))}
{\prod_{j=1}^{n} Df_{i_j} (h_{j-1}(y))} \right) \right|\\
&\leq&
\sum_{j=1}^{n} \big| \log Df_{i_j} (h_{j-1}(x)) - \log Df_{i_j} (h_{j-1}(y)) \big|\\
&\leq&
C \sum_{j=1}^{n} \big| h_{j-1}(x) - h_{j-1}(y) \big|^{1/d}\\
&\leq&
C \sum_{j=0}^{n-1} \big| h_{j}(I) \big|^{1/d}.
\end{eqnarray*}
If we denote by $M_n$ the value of the last sum, this implies that
\begin{equation}\label{general-estimate}
Dg^k (x) \leq \exp(C M_n) \cdot D g^k (h_n(x)).
\end{equation}
In order to control the growth of $S_n$, we will use
the first of the two properties provided by the next

\vspace{0.1cm}

\begin{lem} Let $\ell \!: \mathbb{N}^d_0 \to \mathbb{R}$
be a function taking only positive values. Assume that
\begin{equation}\label{hyp}
\sum_{(i_1,\ldots,i_d) \in \mathbb{N}_0^d}
\!\! \ell(i_1,\ldots,i_d) \esp < \infty.
\end{equation}
Then there exists a constant $B > 0$ such that for each $n \in \mathbb{N}$
there is a geodesic path of length $n$ in $\mathbb{N}^d_0$, say
$\{ (i_1 (j), \ldots, i_d (j)) \!: 0 \leq j \leq n\}$, satisfying
\hspace{0.02cm} $i_1(0) = \ldots = i_d(0) = 0$ \hspace{0.02cm} and 
\begin{equation}\label{first}
\sum_{j=0}^{n-1} \ell(i_1(j),\ldots,i_d(j))^{1/d} \esp \esp
\leq \esp \esp B \esp \big(\log(n+1)\big)^{1-1/d}
\end{equation}
and
\begin{equation}\label{second}
\ell(i_1(n), \ldots, i_d(n)) \esp\esp \leq \esp\esp \frac{B}{(n+1)^{d-1}}.
\end{equation}
\end{lem}

\noindent{\bf Proof.} Denote by $L$ the sum in (\ref{hyp}) above. As in \cite{DKN},
we consider the Markov process on $\mathbb{N}_0^d$ with transition probabilities
$$p \big( (i_1,\ldots,i_d) \mapsto (i_1,\ldots, i_{j-1}, 1+i_j, i_{j+1},\dots,i_d) \big)
:= \frac{1+i_j}{i_1 + \ldots + i_d + d}.$$
For this process, the transition probabilities in $n$ steps are equidistributed along
the $n$-sphere $S_n$ for every $n \geq 1$:
$$i_1 + \cdots + i_d = n \hspace{0.18cm} \Longrightarrow \hspace{0.18cm} \mathbb{P}_n
\big( (0,\ldots,0) \mapsto (i_1,\ldots,i_d) \big) = \frac{1}{|S_n|}.$$
Letting $A_d>0$ be such that \esp $|S_n| \geq A_d \esp (n+1)^{d-1}$ \esp 
for all $n \geq 0$, a direct application of H\"older's inequality yields, 
for every $0 < \tau < 1$,
\begin{eqnarray*}
\mathbb{E} \left( \sum_{j=0}^{n-1} \ell(i_1(j),\ldots,i_d(j))^{\tau}\right)
&=&
\sum_{j=0}^{n-1} \mathbb{E} \left( \ell(i_1(j),\ldots,i_d(j))^{\tau}\right)\\
&=&
\sum_{j=0}^{n-1} \frac{1}{|S_j|} \sum_{(i_1,\ldots,i_d) \in S_j}
\ell(i_1,\ldots,i_d)^{\tau}\\
&\leq&
\left( \sum_{j=0}^{n-1} \sum_{(i_1,\ldots,i_d) \in S_j} \ell(i_1,\ldots,i_d) \right)^{\tau}
\left( \sum_{j=0}^{n-1} \sum_{(i_1,\ldots,i_d) \in S_j}
\left( \frac{1}{|S_j|} \right)^{\frac{1}{1-\tau}} \right)^{1-\tau}\\
&\leq&
L^{\tau} \left( \sum_{j=0}^{n-1} |S_j|
\left( \frac{1}{|S_j|} \right)^{\frac{1}{1-\tau}} \right)^{1-\tau}\\
&=&
L^{\tau} \left( \sum_{j=0}^{n-1}
\left( \frac{1}{|S_j|} \right)^{\frac{\tau}{1-\tau}} \right)^{1-\tau}\\
&\leq&
\frac{L^{\tau}}{A_d^{\tau}} \left( \sum_{j=0}^{n-1}
\frac{1}{(j+1)^{(d-1)\frac{\tau}{1-\tau}}} \right)^{1-\tau}.
\end{eqnarray*}
Now, for $\tau = \frac{1}{d}$, we have \esp $(d-1)\frac{\tau}{1-\tau} = 1$, \esp hence
$$\mathbb{E} \left( \sum_{j=0}^{n-1} \ell(i_1(j),\ldots,i_d(j))^{1/d} \right)
\leq \frac{L^{1/d}}{A_d^{1/d}} \left( \sum_{j=1}^{n} \frac{1}{j} \right)^{1-1/d} \!
\leq \esp \esp \frac{L^{1/d}}{A_d^{1/d}} \esp \big( \log(n+1) \big)^{1-1/d}.$$
A direct application of Chebyshev's inequality then shows
that with probability larger than $2/3$, we must have
$$\sum_{j=0}^{n-1} \ell(i_1(j),\ldots,i_d(j))^{1/d}
\leq \frac{3 \esp L^{1/d}}{A_d^{1/d}} \big( \log(n+1) \big)^{1-1/d}.$$
Moreover, since \esp 
$\sum_{(i_1,\ldots,i_d) \in S_n} \ell (i_1,\ldots,i_d) \leq L$ 
\esp and the arrival probabilities in $n$ steps are equidistributed 
along $S_n$, with probability larger than $2/3$ we must have
$$\ell \big( i_1(n), \ldots, i_d(n) \big) \leq \frac{3 \esp L}{|S_n|}
\leq \frac{3L}{A_d (n+1)^{d-1}}.$$
Thus, letting
$$B := \max \!\Big\{ \frac{3L^{1/d}}{A_d^{1/d}}, \frac{3L}{A_d} \Big\},$$
we have that (\ref{first}) and (\ref{second}) simultaneously hold with
probability larger than $1/3$. This ensures the existence of the
desired geodesic path. $\hfill\square$

\vspace{0.53cm}

Coming back to the proof of Theorem A, we let $\ell$ be the function
that associates to $(i_1,\ldots,i_d)$ the length of the convex closure
of \esp $\bigcup_{i_{d+1} \in \mathbb{Z}} I_{i_1,\ldots,i_d,i_{d+1}}$.
(Notice that this interval coincides with
$f_1^{i_1} \circ \cdots \circ f_{d}^{i_d}(I)$.) \esp
Let $h_n := f_{i_n} \circ \cdots \circ f_{i_1}$ be a random composition
for which (\ref{first}) and (\ref{second}) hold, and let $g_n$ be the
restriction of $g$ to $h_n (I)$. We claim that $g_n$ has no hyperbolic
fixed point. Otherwise, since $g$ commutes with each $f_i$, it would
have a sequence of hyperbolic fixed points (with the same derivative)
accumulating at a limit point, which is clearly impossible.

We are hence under the hypothesis of Borichev's theorem, and an
application of (\ref{general}) in (\ref{general-estimate}) yields
\begin{eqnarray*}
\max_{x \in I} D g^k(x)
&\leq& \exp(C M_n) \cdot D g^k_n (h_n(x))\\
&\leq& \exp \left(C B \esp \big( \log(n+1)\big)^{1-1/d} \right)
\esp \exp (3 \esp C_{g_n} |h_n(I)|^{\alpha} k^{1-\alpha})\\
&\leq& \exp \left(C B \esp \big( \log(n+1)\big)^{1-1/d} \right) \esp
\exp\left(\frac{3 \esp C_{g} B^{\alpha} k^{1-\alpha}}{(n+1)^{(d-1)\alpha}} \right)\!.
\end{eqnarray*}
Taking $n = n_k$ so that \esp $k^{1-\alpha} \sim n^{(d-1)\alpha}$, \esp hence \esp
$\log(k) \sim \log(n_k)$, \esp we obtain for a certain constant $A$ (independent of $k$),
$$\max_{x \in I} Dg^k (x) \leq \exp \left(A \esp \big( \log(k)\big)^{1-1/d} \right).$$
However, since the last expression is of order \esp $o(k)$, \esp this turns
out to be impossible because of the Polterovich-Sodin lemma.


\section{Proof of the critical generalized Kopell lemma for
different regularities via a deterministic argument}
\label{kopell-2}

\hspace{0.45cm} The proof of Theorem B consists in a combination of the ideas of \S \ref{uno}
and \cite{KN}. The case $d=2$ is relatively straightforward. Nevertheless, for larger
$d$, we need a slight but nontrivial modification of the concatenation argument of
\cite{KN}. Just for pedagogical reasons, we independently develop the cases $d=2$,
$d=3$, and the general case $d \geq 3$, so that to introduce the necessary new
ideas in a progressive manner. (Although the reader should have no problem in
passing directly from the case $d=2$ to the general case $d \geq 3$.)

For all cases, we will proceed by contradiction. We assume that \esp
$\alpha_1 + \ldots + \alpha_d = 1$, \esp and we let again $\ell$ \esp be the function
that associates to $(i_1,\ldots,i_d)$ the length of the convex closure of \esp
$\bigcup_{i_{d+1} \in \mathbb{Z}} I_{i_1,\ldots,i_d,i_{d+1}}$. Then we
consider parallelepipeds $Q(n)$ in $\mathbb{N}_0^d$ whose $s^{th}$-side has
length of order $2^{n \alpha_s}$. For such a $Q(n)$, we set
\begin{equation}\label{ele-n}
L_n := \sum_{(i_1,\ldots,i_d) \in Q(n)} \ell(i_1,\ldots,i_d).
\end{equation}

Our task consists in showing that, for an appropriately chosen sequence $(Q(n))$ of
finite multiplicity $M$ (that is, such that no point is contained in more than $M$ of
these parallelepipeds), there is a positive constant $B$ for which there exist
(non-necessarily nonempty)
geodesic segments \esp $\gamma_{1}^1,\gamma_1^2,\ldots, \gamma_1^{d_1},\gamma_2^1,
\gamma_2^2, \ldots,\gamma_2^{d_2},\ldots$, \esp with $d_n \leq d$,
satisfying the following properties:
\begin{itemize}

\item For each $n \geq 1$ and $1 \leq k \leq d_n$, the segment $\gamma_n^k$
is contained in $Q(n)$.

\item Each of these segments intersects the next nonempty one in the sequence above.

\item For certain positive constants $\alpha, D$ and each $n \geq 1$, at least one of the
segments $\gamma_n^1,\ldots,\gamma_n^{d_n}$ contains no less than $2^{n \alpha} / D$ points.

\item Each $\gamma_n^{k}$ is an unidirectional path pointing
in a $s$-direction, with $s = s_{n,k}$, and
\begin{equation}\label{main-1}
\sum_{(i_1,\ldots,i_d) \in \gamma_n^k} \ell(i_1,\ldots,i_d)^{\alpha_s}
\leq B \max \{ L_n^{\alpha_s}, L_{n+1}^{\alpha_s} \}.
\end{equation}
\end{itemize}
We next explain how such a sequence of geodesic segments allows proving Theorem B.
The next paragraphs will be devoted to the constructions of the parallelepipeds as
well as the sequences of geodesics segments satisfying the desired properties in
the corresponding cases.

First of all, notice that the concatenation of the geodesic segments along
intersecting points produces an infinite (non-necessarily geodesic) path
$\gamma: \mathbb{N}_0 \to \mathbb{N}_0^d$. We will assume that $\gamma$ starts
at the origin: if this is not the case, adding an initial segment $\gamma_0 \subset Q(1)$
from the origin to the initial point of $\gamma_1^1$, one may apply the same arguments
below modulo slightly changing the constant $B$.

For each $n \geq 1$, we denote by $N = N(n)$ the entry-time of $\gamma$
into $Q(n+1)$. If we denote by $s(m)$ the direction corresponding to the
jump from $\gamma(m)$ to $\gamma(m+1)$, then (\ref{main-1}) combined
with H\"older's inequality yields
$$\sum_{m=0}^{N} \ell(\gamma(m))^{\alpha_{s(m)}}
\leq B \sum_{m=1}^{n+1} \sum_{k=1}^d L_m^{\alpha_{s_{m,k}}}
\leq
B \sum_{k=1}^d \left( 2\sum_{m=1}^{n+1} L_m \right)^{\alpha_{s_{m,k}}} (n+1)^{1- \alpha}
\leq
2^{\alpha'} d \esp B M^{\alpha'} \esp (n+1)^{1-\alpha},$$
where $\alpha := \min \{ \alpha_1, \ldots, \alpha_d \}$ and
$\alpha' := \max \{ \alpha_1, \ldots, \alpha_d \}$. Moreover, the
assumption on the size of $Q(n)$ easily implies the asymptotic equivalence
\esp $n \sim \log(N)$. \esp (This equivalence will be even more transparent
for the explicit choice of $Q(n)$ further on.) As a consequence, there is a
constant $A' > 0$ such that the previous estimate becomes
$$\sum_{m=0}^{N} \ell \big( \gamma(m) \big)^{\alpha_{s(m)}}
\leq A' \big( \log(N) \big)^{1-\alpha}.$$

The path $\gamma$ induces a sequence $(h_n)$ of compositions of maps from \esp
$\{f_1,f_1^{-1},\ldots,f_d,f_d^{-1} \}$ \esp such that if $I$ denotes the convex
closure of \esp $\bigcup_{i_{d+1}\in\mathbb{Z}} I_{0,\ldots,0,i_{d+1}}$, \esp then
\begin{equation}\label{main-est}
\sum_{m=0}^N \big|h_m (I) \big|^{\alpha_{s(m)}}
\leq A' \big( \log(N) \big)^{1-\alpha}.
\end{equation}
Let us write \esp $h_m = f_{i_m} \circ \cdots \circ f_{i_1}$, \esp where each
$f_{i_j}$ lies in \esp $\{f_1,f_1^{-1},\ldots,f_d,f_d^{-1} \}$, \esp and
let $C$ be a common upper bound for the $\alpha_i$-H\"older constants of
$\log(Df_i), \log(Df_i^{-1})$, where $1 \leq i \leq d$. Given $n > 1$,
let $N'$ be such that $N(n-1) \leq N' \leq N(n) = N$. For each $x,y$
in $I$, estimate (\ref{main-est}) yields
\begin{eqnarray*}
\left| \log \left( \frac{Dh_{N'} (x)}{Dh_{N'} (y)} \right) \right|
&=&
\left| \log \left( \frac{\prod_{m=1}^{N'} Df_{i_m} (h_{m-1}(x))}
{\prod_{m=1}^{N'} Df_{i_m} (h_{m-1}(y))} \right) \right|\\
&\leq&
\sum_{m=1}^{N'} \big| \log Df_{i_m} (h_{m-1}(x)) - \log Df_{i_m} (h_{m-1}(y)) \big|\\
&\leq&
C \sum_{m=1}^{N'} \big| h_{m-1}(x) - h_{m-1}(y) \big|^{\alpha_{s(m-1)}}\\
&\leq&
C \sum_{m=0}^{N} \big| h_{m}(I) \big|^{\alpha_{s(m)}}\\
&\leq&
C A' \big( \log(N) \big)^{1-\alpha}.
\end{eqnarray*}
Moreover, by the third property of our sequence, we may choose $k$ such that
$\gamma_n^k$ contains at least $2^{n\alpha}/D$ points. Since the sum of the values of
$\ell$ along these points is $\leq 1$, such a segment must contain a point at which
the value of $\ell$ is \esp $\leq D/ 2^{n\alpha}$. In other words, we may choose
$N'$ above so that \esp $\big| h_{N'}(I) \big| \leq D / 2^{n \alpha}$.

Denote by $g_{N'}$ the restriction of $g \!:=\! f_{d+1}$ to $h_{N'} (I)$. As in
\S \ref{uno}, the map $g$ (hence $g_{N'}$) cannot have hyperbolic fixed points.
Therefore, taking derivatives in the equality \esp
$g^k \esp = \esp h_{N'}^{-1} \circ g^k \circ h_{N'}$
\esp and using the previous estimate for $y := g^k(x) \in I$,
Borichev's theorem yields
\begin{eqnarray*}
Dg^k (x)
&=&
\frac{Dh_{N'} (x)}{Dh_{N'} (h_{N'}^{-1} g^k
\esp h_{N'} (x))} \cdot D g^k (h_{N'} (x))\\
&=&
\esp \frac{Dh_{N'} (x)}{Dh_{N'} (g^k(x))} \cdot D g^k_{N'} (h_{N'} (x))\\
&\leq&
\exp \left(C A' \esp \big( \log(N)\big)^{1-\alpha} \right) \esp
\exp \left(3 \esp C_{g} |h_{N'}(I)|^{\alpha_{d+1}} k^{1-\alpha_{d+1}} \right)\\
&\leq&
\exp \left(C A' \esp \big( \log(N)\big)^{1-\alpha} \right) \esp
\exp\left(\frac{3 \esp C_{g} D^{\alpha} k^{1-\alpha_{d+1}}}
{2^{n\alpha\alpha_{d+1}}} \right)\!.
\end{eqnarray*}
Take $n = n_k$ so that \esp $k^{1-\alpha_{d+1}} \sim 2^{n\alpha\alpha_{d+1}}$,
\esp hence \esp $n_k \sim \log(k)$. \esp Using the fact that \esp
$\log(N(n_k)) \sim n_k$, \esp we obtain for a certain constant
$A$ (independent of $k$),
$$\max_{x \in I} Dg^k (x)
\leq
\exp \left(A \esp \big( \log(k)\big)^{1-\alpha} \right).$$
However, the last expression is of order \esp $o(k)$, \esp
which is impossible by the Polterovich-Sodin lemma.


\subsection{The case $d=2$}

\hspace{0.45cm} Following \cite{KN}, we let $Q(n)$ be the rectangle
defined as \esp $Q(\!2n+1\!) := [[4^{n\alpha_1} \!, 4^{(n+1)\alpha_1} ]]
\times [[4^{n\alpha_2} \!, 4^{(n+2)\alpha_2}]]$ \esp and \esp
$Q(\!2n+2\!) := [[4^{n\alpha_1}, 4^{(n+2)\alpha_1}]] \times
[[4^{(n+1)\alpha_2}, 4^{(n+2)\alpha_2}]]$, \esp where $[[x,y]]$
stands for the set of integers between $x$ and $y$. Notice that
the multiplicity of the sequence $(Q(n))$ is $4$.

A set of the form $Q(n) \cap \{j \!=\! const \}$ (resp. $Q(n) \cap \{i \!=\! const \}$)
is said to be an {\em horizontal} (resp. {\em vertical}) segment in $Q(n)$. Notice
that the cardinality of this set $\mathcal{H}_n$ (resp. $\mathcal{V}_n$) of
horizontal (resp. vertical) segments is \esp $\geq 2^{n \alpha_2} / D_1$
(resp. \esp $\geq 2^{n \alpha_1} / D_1$) for a certain constant $D_1 > 0$
(independent of $n$). Moreover, there exists a positive constant $D_2$
such that the number of points in each of these horizontal (resp. vertical)
segments is \esp $\geq D_2 2^{n \alpha_1}$ \esp (resp. $\geq D_2 2^{n \alpha_2}$).

Say that an horizontal segment $\gamma$ in $Q(2n+2)$ is {\em good} if
$$\sum_{(i,j) \in \gamma} \ell (i,j) \leq \frac{L_{2n+2}}{|\mathcal{H}_{2n+2}|}.$$
Clearly, there must be at least one good horizontal segment. For such a
segment $\gamma = \gamma_{2n+2}$, H\"older's inequality yields
$$\sum_{(i,j) \in \gamma} \! \ell(i,j)^{\alpha_1} \leq
\left( \frac{L_{2n+2}}{|\mathcal{H}_{2n+2}|} \right)^{\alpha_1} |\gamma|^{1-\alpha_1}
\leq L_{2n+2}^{\alpha_1} \left(\frac{D_1}{2^{(2n+2)\alpha_2}}\right)^{\alpha_1}
(D_2 2^{(2n+2)\alpha_1})^{1-\alpha_1} = D_1^{\alpha_1} D_2^{\alpha_2} L_{2n+2}^{\alpha_1}.$$
Similarly, say that a vertical segment $\gamma$ in $Q(n)$ is good if
$$\sum_{(i,j) \in \gamma} \ell (i,j) \leq \frac{L_n}{|\mathcal{V}_n|}.$$
Again, there must exist a good vertical segment
$\gamma = \gamma_{2n+1} \subset Q(2n+1)$, and for this segment we have
$$\sum_{(i,j) \in \gamma} \! \ell(i,j)^{\alpha_2} \leq
\left( \frac{L_{2n+1}}{|\mathcal{V}_{2n+1}|} \right)^{\alpha_2} |\gamma|^{1-\alpha_2}
\leq L_{2n+1}^{\alpha_2} \left(\frac{D_1}{2^{(2n+1)\alpha_1}}\right)^{\alpha_2}
(D_2 2^{(2n+1)\alpha_2})^{1-\alpha_2} = D_1^{\alpha_2} D_2^{\alpha_1}
L_{2n+1}^{\alpha_2}.$$
Thus, the segments $\gamma_1,\gamma_2,\ldots,\gamma_n, \ldots$ satisfy
(\ref{main-1}) for \esp $B \!\geq\! \max \{D^{\alpha_1}
D_2^{\alpha_2}, D_1^{\alpha_2} D_2^{\alpha_1} \}$. \esp
Each of these segments intersects the next one, and it is easy to check that
between the concatenating points, $\gamma_n$ contains at least $2^{n\alpha} /D$ points
for a certain constant $D > 0$, where $\alpha := \min \{ \alpha_1, \alpha_2 \}$.
Therefore, all the conditions from (\ref{kopell-2}) are fulfilled, and this
concludes the proof of Theorem B in the case $d=2$.


\subsection{The case $d=3$}
\label{d=3}

\hspace{0.45cm} In the case $d \geq 3$, we let $Q(n) := \prod_{i=1}^{d} [[x_n^i,y_n^i]]$
be inductively defined by \esp $Q(1) := [[1,4^d]]^d$ \esp and \esp
$Q(n+1) := \cdots \times [[ 1 + 2^{d\alpha_{m}}(x_n^{m}-1), y_n^m]] \times
[[x_n^{m+1}, 1 + 2^{d\alpha_{m+1}} (y_n^{m+1}-1) ]] \times \cdots $, where for each
$n \geq 1$, we denote by $m = m(n) \in \{1,\ldots, d\}$ the residue class of $n$
modulo $d$. (The dots mean that the corresponding factors remain untouched.)
Notice that the sequence $(Q(n))$ has multiplicity \esp $d+2$.

One easily checks the asymptotic equivalence \esp $y_n^k - x_n^k \sim 2^{n \alpha_k}$.
\esp We let $D_1$ be a constant such that \esp $2^{n \alpha_i} / D_1 \leq
y_n^k - x_n^k \leq D_1 2^{n \alpha_i} - 1$.
\esp Moreover, we may fix a constant $D_2 > 0$ such that \esp
\begin{equation}\label{inocent}
y_{n+1}^{m} - x_{n+1}^m \geq D_2 (y_{n}^{m} - x_{n}^m)
\qquad \mbox{ and } \qquad
y_{n}^{m+1} - x_{n}^{m+1} \geq D_2 (y_{n+1}^{m+1} - x_{n+1}^{m+1}).
\end{equation}

We now specialize to the case $d = 3$. A plane $P$ of the form
\esp $P = Q(n) \cap \{ i_{m+2} = const\}$ \esp will be said to be an
$h$-plane, and the family of $h$-planes in $Q(n)$ will be denoted by $\mathcal{P}_n$.
Notice that the cardinality of $\mathcal{P}_n$ is \esp $\geq 2^{n\alpha_{m+2}} / D_1$.

Given an $h$-plane \esp $P \in \mathcal{P}_{n}$, \esp an {\em horizontal} (resp.
{\em vertical}) segment in $P$ is a set of the form \esp $P \cap \{ i_{m+1} = const \}$
\esp (resp. \esp $P \cap \{i_m = const \}$). The cardinality of the family $\mathcal{H}_n$
(resp. $\mathcal{V}_n$) of horizontal (resp. vertical) segments is \esp
$\geq 2^{n \alpha_{m+1}} / A$ (resp. \esp $\geq 2^{n \alpha_m} / A'$)
for certain positive constants $A,A'$.

Let $\lambda \geq 1$ be fixed. Say that an $h$-plane
$P$ in $\mathcal{P}_n$ is {\em $\lambda$-good} if (see (\ref{ele-n}))
$$\sum_{(i_1,i_2,i_3) \in P} \ell (i_1,i_2,i_3) \leq \frac{\lambda L_n}{|\mathcal{P}_n|}.$$
By Chebyshev's inequality, the proportion of $h$-planes that are $\lambda$-good is
larger than $1 - 1/\lambda$. Similarly, say that an horizontal segment $\gamma$ of
$P \in \mathcal{P}_n$ is {\em $\lambda$-good relatively to $P$} if
$$\sum_{(i_1,i_2,i_3) \in \gamma} \ell (i_1,i_2,i_3) \leq
\frac{\lambda \sum_{(i_1,i_2,i_3) \in P} \esp\esp
\ell(i_1,i_2,i_3) }{|\mathcal{H}_n|}.$$
As before, horizontal directions are relatively $\lambda$-good in a proportion
larger than $1 - 1/\lambda$. Finally, say that a vertical segment
$\gamma$ in $P \in \mathcal{P}_{n}$ is $\lambda$-good if
$$\sum_{(i_1,i_2,i_3) \in \gamma} \ell (i_1,i_2,i_3) \leq
\frac{\lambda \sum_{(i_1,i_2,i_3) \in P} \esp\esp
\ell(i_1,i_2,i_3) }{|\mathcal{V}_n|}.$$
Once again, vertical segments are $\lambda$-good
in a proportion larger than $1 - 1/\lambda$.

Now fix $\lambda \geq 2/D_2$, and let $P \in \mathcal{P}_n$ be a $\lambda$-good $h$-plane. By
the first inequality in (\ref{inocent}), more than a half of the vertical segments of $P$
contained in $Q(n+1)$ are $\lambda$-good. Since $\lambda \geq 2$, more than a half of the
$h$-planes in $Q(n+1)$ are $\lambda$-good. Therefore, there must exist a $\lambda$-good
$h$-plane $P' \in \mathcal{P}_{n+1}$ such that $P \cap P'$ is a $\lambda$-good vertical
segment of $P$. Moreover, by the second inequality in (\ref{inocent}), more than a half
of the vertical segments of $P'$ contained in $Q(n)$ are $\lambda$-good relatively to $P'$.

We may hence fix a sequence $(P_n)$ of $\lambda$-good $h$-planes such that $P_n \cap
P_{n+1}$ is a $\lambda$-good vertical segment $\gamma_n^2$ of $P_n$, for each $n \geq 1$.
(See Figure 1.) Each $P_n$ must contain a relatively $1$-good horizontal segment $\gamma_n^1$.
Finally, let $\gamma_{n}^3$ be a $\lambda$-good vertical segment of $P_{n+1}$ contained in $Q(n)$.
(See Figure 2.) We have thus constructed an infinite sequence of geodesic segments \esp
$\gamma_1^1,\gamma_1^2,\gamma_1^3,\gamma_2^1,\gamma_2^2,\gamma_2^3,\ldots$, \esp each of
which intersects the next one in the sequence. Moreover, since $P_n$ is $\lambda$-good
and $\gamma_n^1$ is a relatively $1$-good horizontal segment in $P_n$,
$$\sum_{(i_1,i_2,i_3) \in \gamma_n^1} \ell (i_1,i_2,i_3)
\leq
\frac{\sum_{(i_1,i_2,i_3) \in P_n} \esp\esp \ell(i_1,i_2,i_3) }{|\mathcal{H}_n|}
\leq
\frac{1}{|\mathcal{H}_n|} \frac{\lambda L_n}{|\mathcal{P}_n|}
\leq \frac{A D_1 \lambda L_n}{2^{n\alpha_{m+1}} 2^{n \alpha_{m+2}}}.$$
By H\"older's inequality, this implies that
$$\sum_{(i_1,i_2,i_3) \in \gamma_n^1} \ell (i_1,i_2,i_3)^{\alpha_{m}}
\leq
\left( \frac{A D_1 \lambda L_n}{2^{n\alpha_{m+1}} 2^{n \alpha_{m+2}}} \right)^{\alpha_{m}}
|\gamma_n^2|^{1 - \alpha_{m}}
\leq
\left( \frac{A D_1 \lambda L_n}{2^{n\alpha_{m+1}} 2^{n \alpha_{m+2}}} \right)^{\alpha_{m}}
(D_1 2^{n\alpha_{m}})^{1 - \alpha_{m}},$$
hence
\begin{equation}\label{first-estim}
\sum_{(i_1,i_2,i_3) \in \gamma_n^1} \ell (i_1,i_2,i_3)^{\alpha_{m}}
\leq (A D_1 \lambda)^{\alpha_{m}} D_1^{1-\alpha_{m}} L_n^{\alpha_{m}}.
\end{equation}
Similarly, since $P_n$ is $\lambda$-good and $\gamma_n^2$ is a $\lambda$-good
vertical segment of $P_n$,
$$\sum_{(i_1,i_2,i_3) \in \gamma_n^2} \ell (i_1,i_2,i_3)
\leq
\frac{\lambda \sum_{(i_1,i_2,i_3) \in P_n} \esp\esp \ell(i_1,i_2,i_3) }{|\mathcal{V}_n|}
\leq
\frac{\lambda}{|\mathcal{V}_n|} \frac{\lambda L_n}{|\mathcal{P}_n|}
\leq \frac{A' D_1 \lambda^2 L_n}{2^{n\alpha_{m}} 2^{n \alpha_{m+2}}}.$$
Again, by H\"older's inequality,
$$\sum_{(i_1,i_2,i_3) \in \gamma_n^2} \!\!\! \ell (i_1,i_2,i_3)^{\alpha_{m+1}}
\leq
\left( \frac{A' D_1 \lambda^2 L_n}{2^{n\alpha_{m}} 2^{n \alpha_{m+2}}} \right)^{\!\alpha_{m+1}}
\! |\gamma_n^2|^{1 - \alpha_{m+1}}
\leq
\left( \frac{A' D_1 \lambda^2 L_n}{2^{n\alpha_{m}} 2^{n \alpha_{m+2}}} \right)^{\!\alpha_{m+1}}
\!\! (D_1 2^{n\alpha_{m+1}})^{1 - \alpha_{m+1}},$$
hence
\begin{equation}\label{second-estim}
\sum_{(i_1,i_2,i_3) \in \gamma_n^k} \ell (i_1,i_2,i_3)^{\alpha_{m+1}}
\leq (A' D_1 \lambda^2)^{\alpha_{m+1}} D_1^{1-\alpha_{m+1}} L_n^{\alpha_{m+1}}.
\end{equation}
Finally, $\gamma_n^3$ is a $\lambda$-good vertical segment of $P_{n+1}$,
which is a $\lambda$-good horizontal plane of $Q(n+1)$, hence
$$\sum_{(i_1,i_2,i_3) \in \gamma_n^3} \ell (i_1,i_2,i_3)
\leq
\frac{\lambda \sum_{(i_1,i_2,i_3) \in P_{n+1}} \esp\esp \ell(i_1,i_2,i_3) }{|\mathcal{V}_{n+1}|}
\leq
\frac{\lambda}{|\mathcal{V}_{n+1}|} \frac{\lambda L_{n+1}}{|\mathcal{P}_{n+1}|}
\leq \frac{A' D_1 \lambda^2 L_{n+1}}{2^{(n+1)\alpha_{m+1}} 2^{(n+1) \alpha_{m}}},$$
and H\"older's inequality yields
\begin{eqnarray*}
\sum_{(i_1,i_2,i_3) \in \gamma_n^3} \!\! \ell (i_1,i_2,i_3)^{\alpha_{m+2}}
&\leq&
\left( \frac{A' D_1 \lambda^2 L_{n+1}}{2^{(n+1)\alpha_{m+1}} 2^{(n+1) \alpha_{m}}}
\right)^{\!\alpha_{m+2}} \! |\gamma_n^3|^{1 - \alpha_{m+2}}\\
&\leq&
\left(\frac{A' D_1 \lambda^2 L_{n+1}}{2^{(n+1)\alpha_{m+1}} 2^{(n+1) \alpha_{m}}}
\right)^{\!\alpha_{m+2}} \! (D_1 2^{(n+1)\alpha_{m+2}})^{1 - \alpha_{m+2}},
\end{eqnarray*}
that is
\begin{equation}\label{third-estim}
\sum_{(i_1,i_2,i_3) \in \gamma_n^3} \ell (i_1,i_2,i_3)^{\alpha_{m+2}} \leq
(A' D_1 \lambda^2)^{\alpha_{m+2}} D_1^{1-\alpha_{m+2}} L_{n+1}^{\alpha_{m+2}}.
\end{equation}

By (\ref{first-estim}), (\ref{second-estim}) and
(\ref{third-estim}), condition (\ref{main-1}) holds for \esp
$B \geq \max_{k} \max \{ (A D_1 \lambda)^{\alpha_{k}} \! D_1^{1-\alpha_k} \!,
(A' D_1 \lambda^2)^{\alpha_k} \! D_1^{1-\alpha_{k}} \}$. \esp Finally, it is
easy to see that for a certain constant $D > 0$, each $\gamma_n^1$ contains
at least $2^{n\alpha}/D$ points between the concatenating points,
where $\alpha := \min \{ \alpha_1, \alpha_2, \alpha_3\}$. This
concludes the checking of the properties from \S \ref{kopell-2}.

\vspace{0.55cm}


\beginpicture

\setcoordinatesystem units <0.59cm,0.59cm>

\putrule from 0 0 to 8 0
\putrule from 0 3 to 8 3
\putrule from 0 0 to 0 3
\putrule from 8 0 to 8 3
\putrule from 10 2 to 10 5
\putrule from 2 5 to 10 5
\putrule from 14 6 to 14 9
\putrule from 5.5 0 to 5.5 3
\putrule from 11.5 9 to 14 9

\putrule from 13 0 to 15.5 0
\plot 15.5 0 15.2 0.15 /
\plot 15.5 0 15.2 -0.15 /

\putrule from 13 0 to 13 2.5
\plot 13 2.5 12.85 2.2 /
\plot 13 2.5 13.15 2.2 /

\plot 13 0 14.5 1.5 /
\plot 14.5 1.5 14.4 1.15 /
\plot 14.5 1.5 14.1 1.34 /

\plot 28.2 6 28.2 9 /
\plot 22.2 0 22.2 3 /
\plot 22.2 0 23 0.8 /
\plot 25 2.8 28.2 6 /
\plot 22.2 3 28.2 9 /
\plot 17 0.8 22.2 0.8 /
\plot 22.2 0.8 25 0.8 /
\plot 19 2.8 22.2 2.8 /
\plot 24.2 2.8 27 2.8 /
\plot 17 0.8 19 2.8 /
\plot 25 0.8 27 2.8 /

\plot 22.2 0.8 24.2 2.8 /
\plot 22.2 0.815 24.2 2.815 /
\plot 22.2 0.83 24.2 2.83 /
\plot 22.2 0.77 24.2 2.77 /
\plot 22.2 0.785 24.2 2.785 /

\plot 17.5 1.3 22.2 1.3 /
\plot 17.5 1.31 22.2 1.31 /
\plot 17.5 1.32 22.2 1.32 /
\plot 17.5 1.29 22.2 1.29 /
\plot 17.5 1.28 22.2 1.28 /

\plot 22.7 1.3 25.5 1.3 /
\plot 22.7 1.31 25.5 1.31 /
\plot 22.7 1.32 25.5 1.32 /
\plot 22.7 1.29 25.5 1.29 /
\plot 22.7 1.28 25.5 1.28 /

\plot 23.8 2.4 23.8 4.6 /
\plot 23.81 2.4 23.81 4.6 /
\plot 23.82 2.4 23.82 4.6 /
\plot 23.78 2.4 23.78 4.6 /
\plot 23.79 2.4 23.79 4.6 /

\plot 22.2 2.5 28.2 8.5 /
\plot 22.2 2.51 28.2 8.51 /
\plot 22.2 2.52 28.2 8.52 /
\plot 22.2 2.49 28.2 8.49 /
\plot 22.2 2.48 28.2 8.48 /

\plot 8 0 14 6 /
\plot 8 3 14 9 /
\plot 5.5 3 11.5 9 /
\plot 0 3 2 5 /
\plot 8 3 10 5 /

\plot 6.2 0.8 8.2 2.8 /
\plot 6.2 0.815 8.2 2.815 /
\plot 6.2 0.77 8.2 2.77 /
\plot 6.2 0.83 8.2 2.83 /
\plot 6.2 0.785 8.2 2.785 /

\setdashes

\plot 12.2 6 12.2 9 /
\plot 6.2 0 6.2 3 /
\plot 6.2 0 12.2 6 /
\plot 6.2 3 12.2 9 /

\plot 0 0.8 8 0.8 /
\plot 2 2.8 10 2.8 /
\plot 0 0.8 2 2.8 /
\plot 8 0.8 10 2.8 /

\plot 22.3 1.3 22.5 1.3 /
\plot 23.8 1.6 23.8 2.4 /
%


\setdots

\plot 0 0 2 2 /
\plot 8 0 10 2 /
\plot 2 2 10 2 /
\plot 2 2 2 5 /
\plot 5.5 0 11.5 6 /
\plot 11.5 6 11.5 9 /
\plot 11.5 6 14 6 /

\plot 22.2 2.8 24.2 2.8 /
\plot 23 0.8 25 2.8 /

\putrule from 0 0.8 to 8 0.8
\putrule from 0.1 0.9 to 8.1 0.9
\putrule from 0.2 1 to 8.2 1
\putrule from 0.3 1.1 to 8.3 1.1
\putrule from 0.4 1.2 to 8.4 1.2
\putrule from 0.5 1.3 to 8.5 1.3
\putrule from 0.6 1.4 to 8.6 1.4
\putrule from 0.7 1.5 to 8.7 1.5
\putrule from 0.8 1.6 to 8.8 1.6
\putrule from 0.9 1.7 to 8.9 1.7
\putrule from 1 1.8 to 9 1.8
\putrule from 1.1 1.9 to 9.1 1.9
\putrule from 1.2 2 to 9.2 2
\putrule from 1.3 2.1 to 9.3 2.1
\putrule from 1.4 2.2 to 9.4 2.2
\putrule from 1.5 2.3 to 9.5 2.3
\putrule from 1.6 2.4 to 9.6 2.4
\putrule from 1.7 2.5 to 9.7 2.5
\putrule from 1.8 2.6 to 9.8 2.6
\putrule from 1.9 2.7 to 9.9 2.7
\putrule from 2 2.8 to 10 2.8

\putrule from 6.2 0 to 6.2 3
\putrule from 6.3 0.1 to 6.3 3.1
\putrule from 6.4 0.2 to 6.4 3.2
\putrule from 6.5 0.3 to 6.5 3.3
\putrule from 6.6 0.4 to 6.6 3.4
\putrule from 6.7 0.5 to 6.7 3.5
\putrule from 6.8 0.6 to 6.8 3.6
\putrule from 6.9 0.7 to 6.9 3.7
\putrule from 7 0.8 to 7 3.8
\putrule from 7.1 0.9 to 7,1 3.9
\putrule from 7.2 1 to 7.2 4
\putrule from 7.3 1.1 to 7.3 4.1
\putrule from 7.4 1.2 to 7.4 4.2
\putrule from 7.5 1.3 to 7.5 4.3
\putrule from 7.6 1.4 to 7.6 4.4
\putrule from 7.7 1.5 to 7.7 4.5
\putrule from 7.8 1.6 to 7.8 4.6
\putrule from 7.9 1.7 to 7.9 4.7
\putrule from 8 1.8 to 8 4.8
\putrule from 8.1 1.9 to 8.1 4.9
\putrule from 8.2 2 to 8.2 5
\putrule from 8.3 2.1 to 8.3 5.1
\putrule from 8.4 2.2 to 8.4 5.2
\putrule from 8.5 2.3 to 8.5 5.3
\putrule from 8.6 2.4 to 8.6 5.4
\putrule from 8.7 2.5 to 8.7 5.5
\putrule from 8.8 2.6 to 8.8 5.6
\putrule from 8.9 2.7 to 8.9 5.7
\putrule from 9 2.8 to 9 5.8
\putrule from 9.1 2.9 to 9.1 5.9
\putrule from 9.2 3 to 9.2 6
\putrule from 9.3 3.1 to 9.3 6.1
\putrule from 9.4 3.2 to 9.4 6.2
\putrule from 9.5 3.3 to 9.5 6.3
\putrule from 9.6 3.4 to 9.6 6.4
\putrule from 9.7 3.5 to 9.7 6.5
\putrule from 9.8 3.6 to 9.8 6.6
\putrule from 9.9 3.7 to 9.9 6.7
\putrule from 10 3.8 to 10 6.8
\putrule from 10.1 3.9 to 10.1 6.9
\putrule from 10.2 4 to 10.2 7
\putrule from 10.3 4.1 to 10.3 7.1
\putrule from 10.4 4.2 to 10.4 7.2
\putrule from 10.5 4.3 to 10.5 7.3
\putrule from 10.6 4.4 to 10.6 7.4
\putrule from 10.7 4.5 to 10.7 7.5
\putrule from 10.8 4.6 to 10.8 7.6
\putrule from 10.9 4.7 to 10.9 7.7
\putrule from 11 4.8 to 11 7.8
\putrule from 11.1 4.9 to 11.1 7.9
\putrule from 11.2 5 to 11.2 8
\putrule from 11.3 5.1 to 11.3 8.1
\putrule from 11.4 5.2 to 11.4 8.2
\putrule from 11.5 5.3 to 11.5 8.3
\putrule from 11.6 5.4 to 11.6 8.4
\putrule from 11.7 5.5 to 11.7 8.5
\putrule from 11.8 5.6 to 11.8 8.6
\putrule from 11.9 5.7 to 11.9 8.7
\putrule from 12 5.8 to 12 8.8
\putrule from 12.1 5.9 to 12.1 8.9
\putrule from 12.2 6 to 12.2 9

\put {$Q(n)$} at  4 5.4
\put{$Q(n+1)$} at 8 7.2
\put {$P_n$} at  0.7 2.3
\put{$P_{n+1}$} at 13.05 6.5

\put{$\gamma_n^1$} at 20 1.85
\put{$\gamma_n^2$} at 22.8 2.2
\put{$\gamma_{n}^3$} at 24.23 3.6
\put{$\gamma_{n+1}^1$} at 27 6.3

\put{Figure 1} at 6.7 -1
\put{Figure 2} at 23.7 -1

\tiny

\put{$(m+2)$-direction} at 13.6 2.9
\put{$(m+1)$-direction} at 15.3 1.85
\put{$m$-direction} at 14.2 -0.4


\endpicture


\subsection{The general case $d \geq 3$}

\hspace{0.45cm} For the proof in the general ($d \geq 3$) case, it is
better to isolate and inductively extend a concatenation argument that
has already been used (in a weak form) in the cases $d = 2$ and $d=3$.

\vspace{0.3cm}

\noindent{\bf The $d$-dimensional black box.} Let $Q$ be a $d$-dimensional
parallelepiped, with $d \geq 3$, and let $\ell$ be a positive function
defined on $Q$. Set
$$L_Q := \sum_{(i_1,\ldots,i_d) \in Q} \ell(i_1,\ldots,i_d).$$
For $1 \leq m \leq d$, an {\em $m$-segment} is a set of the form \esp
$Q \cap \{i_{j} \!=\! const_j \mbox{ except for } j = m \}$. Denote by
$\mathcal{S}(m)$ the family of all $m$-segments. The elements
in the family $\mathcal{S} := \bigcup_m \mathcal{S}(m)$ will be said
to be {\em unidirectional segments}. Given $\lambda \geq 1$,
say that that $\gamma \in \mathcal{S}(m)$ is $\lambda$-good if
$$\sum_{(i_1,\ldots,i_d) \in \gamma} \ell (i_1,\ldots,i_d)
\leq \frac{\lambda L_{Q}}{|\mathcal{S}(m)|}.$$

More generally, let us consider a $d'$-dimensional parallelepiped $Q'$
and a $d''$-dimensional parallelepiped $Q''$, both contained in $Q$,
such that \esp $d' < d''$ \esp and $Q' \subset Q''$. Denote by
$c(Q',Q'')$ the number of (disjoint) translates of $Q'$ that fill
$Q''$ (that is, the number of copies of $Q'$ in $Q''$). For example,
$c(\gamma,Q) = |\mathcal{S}(m)|$ for any $m$-segment $\gamma$. Say
that $Q'$ is $\lambda$-good relatively to $Q''$ whenever
$$\sum_{(i_1,\ldots,i_d) \in Q'} \ell(i_1,\ldots,i_d) \esp\esp \leq \esp\esp
\frac{\lambda}{c(Q',Q'')} \sum_{(i_1,\ldots,i_d) \in Q''} \ell(i_1,\ldots,i_d).$$
Now, say that $\gamma \in \mathcal{S}(m)$ is fully $\lambda$-good if there
exists a {\em flag} of $d'$-dimensional parallelepipeds $Q^{d'}$ that are
$\lambda$-good in $Q$, have the form \esp $Q^{d'} = Q \cap \{ i_{j} = const_j
\mbox{ except for } j \!=\! m,m+1,\ldots,m+d'-1 \}$, \esp and satisfy
$\gamma = Q^1 \subset Q^2 \subset \ldots \subset Q^{d-1}$.

\vspace{0.1cm}

\begin{lem} Given $0 < \kappa < 1$ and $\mu \geq 1$, there exists
$\lambda' = \lambda' (\mu,\kappa,d)$ (not depending on $Q$) such
that the following holds: If $\gamma$ is a fully $\mu$-good
$1$-segment of $Q$, then for a proportion larger than $\kappa$ of
the points $(i_1,\ldots,i_d) \!\in\! Q$, there exists a finite sequence
of $\lambda'$-good unidirectional segments $\gamma_1, \ldots, \gamma_{d'-1}$,
with $d' \leq d$, such that each $\gamma_i$ intersects $\gamma_{i+1}$, with
$\gamma_1$ starting at a point of $\gamma$ and $\gamma_{d'-1}$
ending at $(i_1,\ldots,i_d)$.
\label{bb}
\end{lem}

\noindent{\bf Proof.} We leave the case $d = 3$ to the reader (it
uses similar arguments to those below; compare also \S \ref{d=3}).

Assume inductively that the claim holds in dimension $d$, and
let us deal with the  $(d+1)$-dimensional case. Let \esp
$\gamma = Q^1 \subset \ldots \subset Q^{d}$ \esp be the flag of
$\mu$-good parallelepipeds associated to $\gamma$. Fix $\lambda \geq 1$
large-enough so that \esp $\kappa < (d-1)(1-1/\lambda) - (d-2)$. \esp For
each $2 \leq m \leq d$, Chebyshev's inequality implies that the proportion
of the set of integers $i$ in the projection of $Q$ along the first coordinate
for which $Q^m(i) := Q^m \cap \{ i_1 = i \}$ is $\lambda$-good relatively to
$Q^m$ is larger than $1 - 1/\lambda$. Therefore, for a proportion larger than
\esp $(d-1)(1-1/\lambda) - (d-2)$ \esp of this set of integers $i$, these
properties hold simultaneously, which means that the segment $Q^2 (i)$
is fully $\lambda$-good in $Q^d (i)$. On each such a $Q^d (i)$, the
inductive procedure yields a proportion larger than \esp
$\frac{\kappa }{[(d-1)(1-1/\lambda) - (d-2)]}$ \esp of points in
$Q^d (i)$ that may be reached by concatenating no more than $d$
unidirectional $\lambda'$-good segments $\gamma_2, \ldots, \gamma_{d'}$
of $Q^d (i)$ (with $\gamma_2$ starting at a point of $\gamma_1 := Q^2 (i)$),
where $\lambda' = \lambda' \Big(
\lambda,\frac{\kappa}{[(d-1)(1-1/\lambda) - (d-2)]},d \Big)$.
Notice that each of these segments is \esp $\mu \lambda'$-good
in $Q$. Thus we have a proportion larger than
$$[(d-1)(1-1/\lambda) - (d-2)] \cdot \frac{\kappa}{[(d-1)(1-1/\lambda) - (d-2)]}
= \kappa$$
of points in $Q$ that may be reached by a sequence of $d$ unidirectional
segments that are $\mu \lambda'$-good, the first of which intersects
$\gamma$. This concludes the inductive proof. $\hfill\square$

\vspace{0.43cm}

\noindent{\bf Proof of Theorem B (for $d \geq 3$).} We consider the sequence of
parallelepipeds $Q(n)$ from the begining of \S \ref{d=3}. Fix  $\lambda > 2(d-1)$,
and let $\lambda' := \lambda' (\lambda,1/2,d)$ be the constant defined in the
statement of Lemma \ref{bb}. We will perform a process that starts by arbitrarily
choosing a fully $\lambda$-good $1$-segment $\gamma_1^1$ of $Q(1)$. (Since \esp
$\lambda > 2(d-1)$, \esp we have \esp $(d-1)(1-1/\lambda) - (d-2) > 1/2 > 0$, \esp
and an application of Chebyshev's inequality ensures the existence of such a segment.)

Assume now that there is a concatenating sequence of unidirectional segments
\esp $\gamma_1^1,\ldots,\gamma_1^{d'_1},\ldots,$ $\gamma_{n-1}^{1},
\ldots, \gamma_{n-1}^{d'_{n-1}},\gamma_n^1$, \esp
with each $d'_{j} \leq d$, such that:

\begin{itemize}

\item For $1 \leq n' \leq n-1$ and $2 \leq k \leq d$, the
segment $\gamma_{n'}^k$ is $\lambda'$-good in $Q(n') \cap Q(n'+1)$.

\item For $1 \leq n' \leq n$, the segment $\gamma_{n'}^1$
is a fully $\lambda'$-good $m(n')$-segment in $Q(n')$.

\end{itemize}

We would like to extend this sequence by appropriately choosing
$\gamma_n^2,\ldots,\gamma_n^{d_n'}$ and $\gamma_{n+1}^1$. To do this, we first
invoke Lemma \ref{bb}, which ensures that more than a half of the points of
$Q(n) \cap Q(n+1)$ may be reached starting at a point of $\gamma_n^1$ by
concatenating $\lambda'$-good segments $\gamma_n^2,\ldots,\gamma_n^{d'}$,
with $d' \leq d$.
On the other hand, since \esp $(d-1)(1-1/\lambda) - (d-2) > 1/2 > 0$,
\esp an application of Chebyshev's inequality ensures that for more than
a half of the points of \esp $Q(n) \cap Q(n+1)$ \esp lie in a fully $\lambda$-good
$m(n+1)$-segment $\gamma_{n+1}^1$ of $Q(n+1)$. These two sets must necessarily
intersect, and this allows to define the desired concatenating segments.

Checking the properties from \S \ref{kopell-2} now mimics the cases $d=2$ and $d=3$.
Indeed, let $A > 0$ be a constant such that for every $s$-segment $\gamma$ in $Q$,
$$2^{n (1-\alpha_{s})} / A \leq c(\gamma,Q) \leq A \esp 2^{n (1-\alpha_{s})}.$$
If $\gamma = \gamma_n^k$ and $k \neq 1$, then letting $s:= s_{n,k}$ be
its direction, H\"older's inequality yields
\begin{eqnarray*}
\sum_{(i_1,\ldots,i_d) \in \gamma_n^k} \ell (i_1,\ldots,i_d)^{\alpha_{s}}
&\leq&
\left( \sum_{(i_1,\ldots,i_d) \in \gamma_n^k} \ell (i_1,\ldots,i_d)
\right)^{\!\alpha_{s}} |\gamma_n^k|^{1-\alpha_{s}}\\
&\leq&
\left( \frac{\lambda' L_{Q(n) \cap Q(n+1)}}{c(\gamma_n^k,Q_n)} \right)^{\alpha_{s}}
(A' 2^{n \alpha_{s}})^{1 - \alpha_{s}}\\
&\leq& (\lambda' A)^{\alpha_{s}} (A')^{1-\alpha_{s}} L_n^{\alpha_{s}}.
\end{eqnarray*}
In the case of $\gamma_n^1$, a similar argument applies, so that (\ref{main-1})
holds for \esp
$B \geq \max_k \{ (\lambda'' A)^{\alpha_{k}} (A')^{1-\alpha_{k}} \},$
\esp where \esp $\lambda'' = \max \{ \lambda, \lambda' \}$. Finally, it is easy
to see that for a certain constant $D > 0$, each segment $\gamma_n^1$ contains
at least $2^{n\alpha}/D$ points, where $\alpha :=
\min \{ \alpha_1, \ldots, \alpha_d \}$. $\hfill\square$


\section{Proof of the non-smoothability of the Farb-Franks action in critical regularity}
\label{FF}

\hspace{0.45cm} We now deal with the group $N_d$ of \esp $(d+1) \times (d+1)$ \esp
lower-triangular matrices with integer entries, all of which are equal to 1
on the diagonal. For $i > j$, we denote by $f_{i,j}$ the element represented
by a matrix whose only nonzero entry outside the diagonal is the $(i,j)$-entry
and equals 1. Notice that the $f_{i,j}$'s generate $N_d$.

Let us briefly remind Farb-Franks' action of $N_d$ on $[0,1]$. First, notice that $N_d$
acts linearly on $\mathbb{Z}^{d+1}$ with the affine hyperplane $1 \!\times\! \mathbb{Z}^{d}$
remaining invariant. The thus-induced action on $\mathbb{Z}^{d}$ allows producing an
action on the interval as follows. Let $\big\{I_{i_1,\ldots,i_{d}}\!\!: \esp\esp
(i_1,\ldots,i_{d}) \!\in\! \mathbb{Z}^d \big\}$ be a family of intervals such
that the sum $\sum_{i_1,\ldots,i_{d}} \vert I_{i_1,\ldots,i_{d}}\vert$ is finite,
say equal to 1 after normalization. We join these intervals lexicographically on the
closed interval $[0,1]$, and we identify $f \in N_d$ to the (unique) homeomorphism of
$[0,1]$ that sends affinely the interval $I_{i_1,\ldots,i_d}$ into $I_{f(i_1,\ldots,i_d)}$,
where $f(i_1,\ldots,i_d)$ stands for the action of \esp $f \in N_d$
\esp on \esp $\mathbb{Z}^d \sim \{ 1 \} \times \mathbb{Z}^d$.

As is shown in \cite{CJN}, if we let $\alpha = \alpha(d) := \frac{2}{d(d-1)}$,
then for every $\varepsilon > 0$, this action is conjugated to an action by
$C^{1 + \alpha - \varepsilon}$ diffeomorphisms, but it cannot be (semi-)conjugated
to an action by $C^{1 + \alpha + \varepsilon}$ diffeomorphisms. Our aim is to extend the
last result to the critical regularity $C^{1+\alpha}$. To do this, we will follow a
similar strategy to that of the generalized Kopell's lemma. We assume that a topological
conjugacy exists, and for simplicity we continue denoting by $I_{i_1,\ldots,i_d}$ the
image of the corresponding interval under this conjugacy. We let $I$ be the convex
closure of \esp $\bigcup_{i_d \in \mathbb{Z}} I_{0,\ldots,0,i_d}$. Notice that the
element $g := f_{d+1,1}$ lies in the center of $N_d$ and fixes the interval $I$.
Moreover, every element in $N_{d}$ sends $I$ into either itself or a disjoint interval.
If we consider the isomorphic copy $N_{d-1}^* \subset N_d$ of $N_{d-1}$ formed by all
elements whose last row and column coincide with those of the identity, then the orbit
of $I$ under $N_d$ coincides with that under $N_{d-1}^*$. Moreover, the stabilizer of
$I$ under the $N_{d-1}^*$-action corresponds to the subgroup formed by the elements
whose first column coincides with those of the identity. Since this subgroup is
naturally isomorphic to $N_{d-2}$, the orbit-graph of $I$ identifies to a coset
space $N_{d-1} / N_{d-2}$, and has $\mathbb{Z}^{d-1}$ as set of vertices.
(See \cite[Figure 2]{CJN} for an illustration in the case $d \!=\! 3$.)


\subsection{From sublinear distortion to the proof of Theorem C}
\label{in-general}

\hspace{0.45cm} As in previous sections, we will decompose (part of) the orbit
of $I$ (which identifies to $\mathbb{Z}^{d-1}$) into parallelepipeds. Following
\cite[\S 2.4]{CJN}, we define $Q(n)$ by induction. We first let
$Q(0) := [ 1,1+4^{d+1} ]^{d-1}$. Now, assuming that
\esp $Q(n):=[x_{1,n},y_{1,n}] \times \cdots \times [x_{d-1,n},y_{d-1,n}]$
\esp has been already defined, we let $i(n) \!\in\! \{1,\ldots,d-1\}$ be the
residue class (mod. $d-1$) of $n$, and we set
$$Q(n+1) := \cdots \times [1+4^{i(n)} (x_{i(n),n}-1), y_{i(n),n}]
\times [x_{i(n)+1,n},1 + 4^{i(n)+1} (y_{i(n)+1,n}-1)] \times \cdots,$$
where the dots mean that the corresponding factors remain untouched.
Notice that all $x_{i,n}$, $y_{i,n}$, and $y_{i,n} - x_{i,n}$,
are asymptotically equivalent to \esp $4^{\frac{in}{d-1}}$.

\vspace{0.1cm}

For each $(i_1,\ldots,i_{d-1}) \in \mathbb{Z}^{d-1}$, we let
$\ell(i_1,\ldots,i_{d-1})$ be the length of the interval $I_{i_1,\ldots,i_{d-1}}$ defined
as the convex closure of \esp $\bigcup_{i_d \in \mathbb{Z}} I_{i_1,\ldots,i_{d-1},i_d}$.
\esp We also set
$$L_n \esp\esp :=
\sum_{(i_1,\ldots,i_{d-1}) \in Q(n)} \ell(i_1,\ldots,i_{d-1}).$$
Our task now consists in showing that there exists a sequence of paths
(segments) \esp $\gamma_0,\gamma_1^1,\ldots,\gamma_1^{k_1}, \ldots,$
$\gamma_n^1,\ldots,\gamma_n^{k_n},\ldots$, \esp with each $k_j \leq K_d$
for a certain constant $K_d$, such that

\begin{itemize}

\item Each $\gamma_n^k$ is contained in $Q(n)$,
whereas $\gamma_0$ is contained in $Q(0)$.

\item For each $n,k$, there exists a generator $f_{i,j}$ of $N_{d-1}^*$ such that
two consecutive points in $\gamma_n^k$ differ by the action of either $f_{i,j}$
or its inverse.

\item There exists a constant $D > 0$ such that for each $n$, at least one of
the $\gamma_n^{k}$ has no less than \esp $4^{\frac{n}{d-1}} / D$ \esp points.

\item There exists a constant $B > 0$ such that for all $n,k$,
\begin{equation}\label{it-works}
\sum_{(i_1,\ldots,i_{d-1}) \in \gamma_n^k} \ell (i_1, \ldots, i_{d-1})^{\alpha}
\leq B L_n^{\alpha}.
\end{equation}

\end{itemize}

Assuming this, we next explain how to complete the proof of Theorem C along the lines
of the arguments given for Theorem B. (Showing the existence of the desired sequences
of parallelepipeds and segments will be postponed to the next two sections.)

The concatenation of the segments above produces an infinite path $\gamma \!:
\mathbb{N}_0 \to \mathbb{N}_0^d$, which we may assume to start at the origin.
(Otherwise, we add an extra initial segment and we slightly change the
constant $B$.) For each $m \geq 0$, we let $f_{m}$ be the element of the
form $f_{i,j}^{\pm 1}$ that moves the $m^{th}$ point of $\gamma$ into the
$(m+1)^{th}$ one, and we denote \esp $h_m := f_{m} \circ \cdots \circ f_1$,
with $h_0 := Id$. Moreover, for each $n \geq 1$, we denote by $N \!=\! N(n)$ the
entry-time of $\gamma$ into $Q(n+1)$. Due to the asymptotics of the lengths of
the sides of $Q(n)$, we have \esp $n \sim  \log(N)$. \esp By (\ref{it-works})
and H\"older's inequality, for a certain constant $A' > 0$, we have
$$\sum_{m=0}^N \! \big|h_m (I) \big|^{\alpha} \!
\leq B \! \sum_{m=0}^{n} \sum_{k=1}^{K_d} L_m^{\alpha}
\leq
B \! \sum_{k=1}^{K_d} \left( \sum_{m=0}^{n} L_m \right)^{\!\!\alpha} \!\! (n+1)^{1- \alpha}
\leq
B K_d (d+2)^{\alpha} (n+1)^{1-\alpha}
\leq A' \big( \log(N) \big)^{1-\alpha}\!\!,$$
where the factor $(d+2)$ comes from the multiplicity of the sequence $(Q(n))$.

Now, for every $x \!\in\! I$, the equality \esp
$g^k = h_m^{-1} \circ g^k \circ h_m$ \esp yields
\begin{equation}\label{eq}
Dg^k (x) = \frac{Dh_m (x)}{Dh_m (y)} \cdot Dg^k (h_m (x)),
\end{equation}
where \esp $y := y_k = g^k(x)$. \esp Since $y$ belongs to $I$, letting
$C$ be a common upper bound for the $\alpha_i$-H\"older constants of
$\log(Df_{i,j}), \log(Df_{i,j}^{-1})$, where $i > j$, for each $N'$
such that $N(n-1) \leq N' \leq N(n) = N$, we have
\begin{eqnarray*}
\left| \log \left( \frac{Dh_{N'} (x)}{Dh_{N'} (y)} \right) \right|
&=&
\left| \log \left( \frac{\prod_{m=1}^{N'} Df_{i_m} (h_{m-1}(x))}
{\prod_{m=1}^{N'} Df_{i_m} (h_{m-1}(y))} \right) \right|\\
&\leq&
\sum_{m=1}^{N'} \big| \log Df_{i_m} (h_{m-1}(x)) - \log Df_{i_m} (h_{m-1}(y)) \big|\\
&\leq&
C \sum_{m=1}^{N'} \big| h_{m-1}(x) - h_{m-1}(y) \big|^{\alpha}\\
&\leq&
C \sum_{m=0}^{N} \big| h_{m}(I) \big|^{\alpha}\\
&\leq&
C A' \big( \log(N) \big)^{1-\alpha}.
\end{eqnarray*}
Moreover, since at least one of the $\gamma_{n}^{k}$'s is assumed
to have no less than $4^{\frac{n}{d-1}} / D$ points, we may choose such an
$N'$ so that \esp $\big| h_{N'}(I) \big| \leq D / 4^{\frac{n}{d-1}}$.
\esp Using Borichev's thorem, this yields
$$Dg^k (h_{N'}(x)) \leq \exp \big( 3 C_g \big| h_{N'}(I) \big|^{\alpha} k^{1-\alpha})
\leq \exp \left( \frac{3 D^{\alpha} C_g k^{1-\alpha}}{4^{\frac{n\alpha}{d-1}}} \right),$$
which due to (\ref{eq}) and the previous estimate implies
$$Dg^k (x) \leq \exp (CA' (\log(N))^{1-\alpha})
\exp \left( \frac{3 D^{\alpha} C_g k^{1-\alpha}}{4^{\frac{n\alpha}{d-1}}} \right).$$
Choosing $n = n_k$ such that \esp $k^{1-\alpha} \sim 4^{\frac{n\alpha}{d-1}}$, \esp hence
\esp $n \sim \log(N) \sim \log(k)$, \esp this yields, for a certain constant $A > 0$,
$$\max_{x \in I} Dg^k (x)
\leq
\exp \left(A \esp \big( \log(k)\big)^{1-\alpha} \right).$$
However, the last expression is of order \esp $o(k)$, \esp
which is impossible by the Polterovich-Sodin lemma.


\subsection{The case $d=3$}
\label{FF:d=3}

\hspace{0.45cm} Again for pedagogical reasons, we first deal with the case $d=3$,
though the reader should have no problem in passing directly to the general case
treated in the next section. Notice that for $d\!=\!3$, the critical value of
$\alpha$ is $1/3$. In analogy to \S \ref{d=3}, let us introduce some terminology.

\vspace{0.1cm}

An {\em horizontal set} in $Q(2n+1)$ is a subset $P = P_r$ of the form
$$Q(2n+1) \cap \{(i,j) \!: \esp i \in [x_{1,2n+1}, y_{1,2n+1}], \esp j \in
[x_{2,2n+1} + (r-1) y_{1,2n+1}, x_{2,2n+1} + r y_{1,2n+1} [ \},$$
where $r \!\in\! \{1,2 \ldots, r_{2n+1}\}$, with
$r_{2n+1} \sim (y_{2,2n+1} - x_{2,2n+1}) / y_{1,2n+1}$ chosen as the
smallest possible index so that $Q(2n+1)$ is the union of the $P_r$'s.
\esp Given $\lambda \geq 1$, such a set is said to be
{\em $\lambda$-good} whenever $r < r_{2n+1}$ and
\begin{equation}\label{hor-2n+1-uno}
\sum_{(i,j) \in P} \ell (i,j)
\leq \frac{\lambda L_{2n+1}}{r_{2n+1}}.
\end{equation}
An {\em horizontal segment} in $Q(2n+1)$ is a subset of the form \esp $Q(2n+1) \cap \{(i,j)
\! : j = const \}$. Such a segment $\gamma$ will be said to be {\em $\lambda$-good relatively
to the horizontal set $P$} containing it whenever
\begin{equation}\label{hor-2n+1-dos}
\sum_{(i,j) \in \gamma} \ell (i,j)
\leq
\frac{\lambda}{y_{1,2n+1}} \sum_{(i,j) \in P} \ell (i,j).
\end{equation}

A {\em vertical set in $Q(2n) \cap Q(2n+1)$} is a set of type
\esp $P_{2n}^{2n+1}(k) := Q(2n) \cap Q(2n+1) \cap \{(i,j) \!: i = k \}$.
A {\em vertical segment \esp in $Q(2n) \cap Q(2n+1)$} is a set of the form \esp
$\gamma_{2n}^{2n+1}(k,r) := P_r \cap P_{2n}^{2n+1}(k)$. \esp This segment is $\lambda$-good
relatively to the vertical set $P_{2n}^{2n+1}(k)$ in \esp $Q(2n) \cap Q(2n+1)$ \esp
containing it whenever
\begin{equation}\label{vert-2n}
\sum_{(i,j) \in \gamma_{2n}^{2n+1}(k,r)} \ell (i,j)
\leq
\frac{\lambda}{r_{2n+1}} \sum_{(i,j) \in P_{2n}^{2n+1}(k)} \!\!\!\! \ell (i,j)
\end{equation}

A {\em vertical set in \esp $Q(2n+1) \cap Q(2n+2)$} \esp is a set \esp 
$P_{2n+1}^{2n+2}(k) := Q(2n+1) \cap Q(2n+2) \cap \{(i,j) \!: i = k \}$.
A {\em vertical segment in $Q(2n+1) \cap Q(2n+2)$} is a set of the form
\esp $\gamma_{2n+1}^{2n+2}(k,r) := P_r \cap P_{2n+1}^{2n+2}(k)$.
\esp This segment is $\lambda$-good relatively to the
horizontal set $P_r$ in $Q(2n+1)$ containing it whenever
\begin{equation}\label{two}
\sum_{(i,j) \in \gamma_{2n+1}^{2n+2}(k,r)} \!\!\!\!\!\!\!\ell (i,j)
\leq
\frac{\lambda}{1 + y_{1,2n+1} - x_{1,2n+1}}
\sum_{(i,j) \in P_r} \!\!\! \ell(i,j).
\end{equation}

Finally, a {\em vertical set in $Q(2n+2)$} is a set of type
$P(k) := Q(2n+2) \cap \{(i,j) \!: i = k \}$.
Such a set $P$ is $\lambda$-good provided
\begin{equation}\label{ufff}
\sum_{(i,j) \in P} \ell (i,j)
\leq
\frac{\lambda L_{2n+2}}{1 + y_{1,2n+2} - x_{1,2n+2}}.
\end{equation}

Now, for each $k \in [x_{1,2n+2},y_{1,2n+2}]$, we decompose
$\{k\} \times [x_{2,2n+2},y_{2,2n+2}] \sim [x_{2,2n+2},y_{2,2n+2}]$
into $k$ paths, each of which has consecutive points at distance $k$.
The resulting paths will be said to be {\em vertical segments in
$Q(2n+2)$}.\footnote{Rather surprisingly, there is no need of the intricate
decomposition procedure of \cite[\S 2.5]{CJN} here.} Such a vertical segment $\gamma$ is
said to be $\lambda$-good relatively to the vertical set $P = P(i)$ in $Q(2n+2)$ containing it if
\begin{equation}\label{vertical-segment}
\sum_{(i,j) \in \gamma} \ell (i,j)
\leq
\frac{\lambda}{i} \sum_{(i,j) \in P} \ell (i,j).
\end{equation}
(Notice that vertical segments in $Q(2n+2)$ naturally arise from the action of $f_{3,2}$.)

\vspace{0.1cm}

Assume we are given a $\lambda$-good vertical set
$P = P(k)$ in $Q(2n)$ and a $1$-good vertical segment $\gamma_{2n}^{1}$
relatively to $P$. For at least a half of the $r \in \{1,2,\ldots,r_{2n+1}-1\}$,
the vertical segment $\gamma_{2n}^{2n+1}(k,r)$ is 2-good relatively to $P_{2n}^{2n+1}(k)$.
Similarly, at least a half of the horizontal sets in $Q(2n+1)$ are 2-good. Consequently,
there must be some $r \in \{1,2,\ldots,r_{2n+1}-1\}$ such that the corresponding vertical segment
$\gamma_{2n}^{2n+1}(k,r) \subset Q(2n) \cap Q(2n+1)$ and horizontal set $P_r \subset Q(2n+1)$
are 2-good. Let $\gamma_{2n+1}^1$ be a 1-good horizontal segment in $P_r$. The segments 
$\gamma_{2n}^1$ and $\gamma_{2n+1}^1$ do not necessarily intersect, but using the 
vertical segment $\gamma_{2n}^2 := \gamma_{2n}^{2n+1} (k,r)$, we can concatenate them.

Assume now that we are given $r$ such that $P_r$ is a $\lambda$-good horizontal set in
$Q(2n+1)$ together with a $1$-good horizontal segment $\gamma_{2n+1}^1$ relatively to
$P_r$. For more than a half of the $k' \in [x_{1,2n+2},y_{1,2n+2}]$, the vertical segment
$\gamma_{2n+1}^{2n+2}(k',r)$ is $2$-good relatively to $P_r$. Similarly, for more than a
half of these $k'$, the vertical set $P(k') \subset Q(2n+2)$ is $2$-good. Take $k'$ lying
simultaneously in both sets, and choose any vertical segment $\gamma_{2n+2}^1$ that is
1-good relatively to $P(k')$. Again, the segments $\gamma_{2n+1}^1$ and $\gamma_{2n+2}^1$
do not necessarily intersect, but using $\gamma_{2n+1}^2 := \gamma_{2n+1}^{2n+2}(k',r)$,
we can concatenate them.

\vspace{0.5cm}


\beginpicture

\setcoordinatesystem units <1.2cm,1.2cm>

\circulararc -180 degrees from 1 0.3 center at 1 0.6
\circulararc -180 degrees from 1 0.9 center at 1 1.2
\circulararc -180 degrees from 1 1.5 center at 1 1.8
\circulararc -180 degrees from 1 2.1 center at 1 2.4
\circulararc -180 degrees from 1 2.7 center at 1 3
\circulararc -180 degrees from 1 3.3 center at 1 3.6

\circulararc -180 degrees from 7 2.2 center at 7 2.6
\circulararc -180 degrees from 7 3 center at 7 3.4
\circulararc -180 degrees from 7 3.8 center at 7 4.2
\circulararc -180 degrees from 7 4.6 center at 7 5
\circulararc -180 degrees from 7 5.4 center at 7 5.8
\circulararc -180 degrees from 7 6.2 center at 7 6.6
\circulararc -180 degrees from 7 7 center at 7 7.4
\circulararc -180 degrees from 7 7.8 center at 7 8.2

\putrule from 0 0 to 2 0
\putrule from 0 0 to 0 4
\putrule from 2 0 to 2 4
\putrule from 0 4 to 8 4
\putrule from 0 2 to 8 2
\putrule from 8 2 to 8 9
\putrule from 5 2 to 5 9
\putrule from 5 9 to 8 9

\putrule from 1 0 to 1 4
\putrule from 1.01 0 to 1.01 4
\putrule from 1.015 0 to 1.015 4
\putrule from 0.985 0 to 0.985 4
\putrule from 0.99 0 to 0.99 4

\putrule from 7 2 to 7 9
\putrule from 7.01 2 to 7.01 9
\putrule from 7.015 2 to 7.015 9
\putrule from 6.985 2 to 6.985 9
\putrule from 6.99 2 to 6.99 9

\putrule from 0 3.2 to 8 3.2
\putrule from 0 3.21 to 8 3.21
\putrule from 0 3.22 to 8 3.22
\putrule from 0 3.23 to 8 3.23

\putrule from 0.98 2.88 to 0.98 3.32
\putrule from 0.97 2.88 to 0.97 3.32
\putrule from 0.96 2.88 to 0.96 3.32
\putrule from 0.95 2.88 to 0.95 3.32
\putrule from 1.05 2.88 to 1.05 3.32
\putrule from 1.04 2.88 to 1.04 3.32
\putrule from 1.03 2.88 to 1.03 3.32
\putrule from 1.02 2.88 to 1.02 3.32

\putrule from 6.98 2.88 to 6.98 3.32
\putrule from 6.97 2.88 to 6.97 3.32
\putrule from 6.96 2.88 to 6.96 3.32
\putrule from 6.95 2.88 to 6.95 3.32
\putrule from 7.05 2.88 to 7.05 3.32
\putrule from 7.04 2.88 to 7.04 3.32
\putrule from 7.03 2.88 to 7.03 3.32
\putrule from 7.02 2.88 to 7.02 3.32


\setdots

\putrule from 0 2.9 to 8 2.9
\putrule from 0 3 to 8 3
\putrule from 0 3.1 to 8 3.1
\putrule from 0 3.2 to 8 3.2
\putrule from 0 3.3 to 8 3.3

\putrule from 0.07 2.95 to 8 2.95
\putrule from 0.07 3.05 to 8 3.05
\putrule from 0.07 3.15 to 8 3.15
\putrule from 0.07 3.25 to 8 3.25

\small

\put{Figure 2} at 4 -0.42
\put{$Q(2n)$} at -0.5 1.3
\put{$Q(2n+1)$} at 3 4.2
\put{$Q(2n+2)$} at 8.78 6
\put{$P_r$} at 3.8 2.7
\put{$P(k)$} at 1 -0.25
\put{$P(k')$} at 7 1.72
\put{$\gamma_{2n+1}^1$} at 3.2 3.46
\put{$\gamma_{2n}^1$} at 0.45 0.8
\put{$\gamma_{2n+2}^1$} at 6.2 6.9
\put{$\gamma_{2n}^2$} at 1.3 2.9
\put{$\gamma_{2n+1}^2$} at 7.47 2.88
\put{} at -3 0

\endpicture


\vspace{0.5cm}

Thus, starting with any vertical segment $\gamma_0$ that is 1-good relatively
to a $1$-good vertical set in $Q_0$, we can produce a concatenating sequence
$\gamma_0,\gamma_1^1,\gamma_1^2,\gamma_2^1,\gamma_2^2,\ldots$. We claim that
this induces a good-enough
sequence of segments in that they satisfy the properties of \S \ref{in-general},
thus concluding the proof of Theorem C for the case $d=3$. Indeed, the first two
properties are clear from the construction, whereas the third one is easily
seen to hold for $\gamma_n^1$. To check the fourth property, that is,
(\ref{it-works}), we will use throughout the asymptotics of
$x_{i,n}, y_{i,n}, y_{i,n}-x_{i,n}$ (which are all of order
$4^{\frac{in}{d-1}}$). Recall also that $\alpha = 1/3$.

For $\gamma_{2n}^1$, using the adapted version of (\ref{ufff}) and
(\ref{vertical-segment}) together with H\"older's inequality, we get
\begin{eqnarray*}
\sum_{(i,j) \in \gamma_{2n}^1} \ell(i,j)^{\alpha}
&\leq&
\left( \frac{2}{i} \sum_{(i,j) \in P (k)} \ell(i,j)
\right)^{\alpha} |\gamma_{2n}^1|^{1-\alpha}\\
&\leq&
\left( \frac{2}{i} \cdot \frac{2 L_{2n}}{1+y_{1,2n}-x_{1,2n}} \right)^{\alpha}
C \! \left( \frac{y_{2,2n} - x_{2,2n}}{i} \right)^{1-\alpha}\\
&\leq&
\frac{C' \esp (y_{2,2n}-x_{2,2n})^{1-\alpha}}{x_{1,2n} \esp
(y_{1,2n}-x_{1,2n})^{\alpha}} \cdot L_{2n}^{\alpha}\\
&\leq&
B \frac{4^{2n(1-\alpha)}}{4^{n} 4^{n\alpha}} L_{2n}^{\alpha}\\
&=&
B \frac{4^{4n/3}}{4^{n} 4^{n/3}} L_{2n}^{\alpha}\\
&=&
B L_{2n}^{\alpha}.
\end{eqnarray*}
For $\gamma_{2n}^2$, using (\ref{vert-2n}) and H\"older's inequality, we get
\begin{eqnarray*}
\sum_{(i,j) \in \gamma_{2n}^2} \ell(i,j)^{\alpha}
&\leq&
\left( \frac{2}{r_{2n+1}} \sum_{(i,j) \in P_{2n}^{2n+1} (k)} \ell(i,j)
\right)^{\alpha} |\gamma_{2n}^2|^{1-\alpha}\\
&\leq&
\left( \frac{2 C y_{1,2n+1}}{y_{2,2n+1} - x_{2,2n+1}}
\cdot \frac{2 L_{2n}}{1+y_{1,2n}-x_{1,2n}} \right)^{\alpha}
y_{1,2n+1}^{1-\alpha}\\
&\leq&
C' \frac{4^{n}}{(4^{(2n+1)} 4^{n})^{\alpha}}\cdot L_{2n}^{\alpha}\\
&\leq&
B \frac{4^{n}}{4^{3n\alpha}} L_{2n}^{\alpha}\\
&=&
B L_{2n}^{\alpha}.
\end{eqnarray*}
For $\gamma_{2n+1}^1$, using the appropriate versions of (\ref{hor-2n+1-uno})
and (\ref{hor-2n+1-dos}), H\"older's inequality yields
\begin{eqnarray*}
\sum_{(i,j) \in \gamma_{2n+1}^1} \ell(i,j)^{\alpha}
&\leq&
\left( \frac{2}{y_{1,2n+1}} \sum_{(i,j) \in P_r} \ell(i,j)
\right)^{\alpha} |\gamma_{2n+1}^1|^{1-\alpha}\\
&\leq&
\left( \frac{2}{y_{1,2n+1}}
\cdot \frac{2 L_{2n+1}}{r_{2n+1}} \right)^{\alpha}
C \esp (1 + y_{1,2n+1} - x_{1,2n+1})^{1-\alpha}\\
&\leq&
C' \frac{(1+ y_{1,2n+1} - x_{1,2n+1})^{(1-\alpha)} y_{1,2n+1}^{\alpha}}
{y_{1,2n+1}^{\alpha} (1+ y_{2,2n+1} - x_{2,2n+1})^{\alpha}} \cdot L_{2n+1}^{\alpha}\\
&\leq&
B \frac{4^{n(1-\alpha)}}{4^{(2n+1)\alpha}} L_{2n+1}^{\alpha}\\
&=&
B L_{2n+1}^{\alpha}.
\end{eqnarray*}
Finally, for $\gamma_{2n+1}^2$, using (\ref{two}) we obtain
\begin{eqnarray*}
\sum_{(i,j) \in \gamma_{2n+1}^2} \ell(i,j)^{\alpha}
&\leq&
\left( \frac{2}{1+y_{1,2n+1}-x_{1,2n+1}} \sum_{(i,j) \in P_r} \ell(i,j)
\right)^{\alpha} |\gamma_{2n+1}^2|^{1-\alpha}\\
&\leq&
\left( \frac{C L_{2n+1}}{(1 + y_{1,2n+1} - x_{1,2n+1})r_{2n+1}} \right)^{\alpha}
y_{1,2n+1}^{1-\alpha}\\
&\leq&
C' \frac{y_{1,2n+1}^{\alpha} y_{1,2n+1}^{1 - \alpha}}
{(1 + y_{2,2n+1} - x_{2,2n+1})^{\alpha} (1 + y_{1,2n+1} - x_{1,2n+1})^{\alpha}}
\cdot L_{2n+1}^{\alpha}\\
&\leq&
B \frac{4^{n}}{4^{2n\alpha} 4^{n\alpha}} L_{2n+1}^{\alpha}\\
&=&
B L_{2n+1}^{\alpha}.
\end{eqnarray*}


\subsection{The general case}
\label{FF:general}

\hspace{0.45cm} As in the case of Theorem B,
to prove Theorem C we will argue by induction.

\vspace{0.3cm}

\noindent{\bf The vertical subdivision procedure.} Given $d \!\geq\! 3$,
we let \esp $Q := \prod_{k=1}^{d-1} [x_k,y_k]$ \esp be a parallelepiped in
$\mathbb{Z}^{d-1}$, where all \esp $x_1,y_1, \ldots, x_{d-1}, y_{d-1}$ \esp
are integers. Given $A \geq 1$, we say that $Q$ is {\em $A$-round} if
\begin{equation}\label{efe-de}
\frac{(1+y_{1}-x_{1})^i}{A}
\leq x_{i}
< y_{i}
\leq A (1 + y_{1} - x_{1})^i,
\quad
\frac{(1 + y_{1} - x_{1})^i}{A}
\leq 1 + y_{i} - x_{i}
\leq A (1 + y_{1} - x_{1})^i\!.
\end{equation}

By cutting along the last coordinate, every $A$-round parallelepiped $Q$ may be
divided into disjoint parallelepipeds $Q_1,\ldots,Q_{M_1}$, each of which has
$(d-1)^{th}$-side of length $y_{d-2} - 1$ possibly excepting the last one.
By (\ref{efe-de}), we have
$$\frac{1+y_1-x_1}{A^2} \leq \frac{1+y_{d-1}-x_{d-1}}{y_{d-2}-1} \leq
M_1 \leq 1 + \frac{1 + y_{d-1} - x_{d-1}}{y_{d-2}-1} \leq 1 + A^2 (1+y_1-x_1).$$
Similarly, each $Q_{m_1}$ satisfying $m_1 < M_1$ may be subdivided into disjoint
parallelepipeds $Q_{m_1,1}$, $Q_{m_1,2},\ldots,Q_{m_1,M_2}$, each of which has
$(d-1)^{th}$-side of length $y_{d-3}-1$ possibly excepting the last one. Again,
(\ref{efe-de}) implies that
$$\frac{1+y_1-x_1}{A^2} \leq M_2 \leq 1 + A^2(1+y_1-x_1).$$
In general, for $k \leq d-2$, each parallelepiped $Q_{m_1,\ldots,m_{k-1}}$ satisfying
$m_j \neq M_j$ for all $j \leq k-1$ may be divided into $Q_{m_1,\ldots,m_{k-1},1},
Q_{m_1,\ldots,m_{k-1},2}, \ldots, Q_{m_1,\ldots,m_{k-1},M_k}$, each of which has
$(d-1)^{th}$-side of length $y_{d-k-2}$ possibly excepting the last one. Moreover,
(\ref{efe-de}) implies that
\begin{equation}\label{general-est}
\frac{1+y_1-x_1}{A^2} \leq M_k \leq 1 + A^2(1+y_1-x_1).
\end{equation}
Here, for $k= 0$, we interpret $Q_{m_1,\ldots,m_k}$ as $Q$.

A {\em level} in $Q$ is a set of the form \esp $H_i \cap Q$, \esp where \esp
$H_i := \{ (i_1,\ldots,i_{d-2},i_{d-1}) \in \mathbb{Z}^{d-1} \!: i_{d-1} = i\}$.
\esp To each level there is a unique associated sequence
\begin{equation}\label{asoc-sequence}
H_i \cap Q \subset Q_{m_1,\ldots,m_{d-2}} \subset Q_{m_1,\ldots,m_{d-3}}
\subset \ldots \subset Q_{m_1} \subset Q.
\end{equation}
Say that the level is {\em admissible} if each of the $m_i'$s above
differs from the corresponding $M_i$. Using (\ref{efe-de}) and (\ref{general-est}),
one easily checks that for a certain constant $A' = A' (A,d) $, the proportion of
non-admissible levels is no larger than
$$\frac{1}{(1+y_{d-1}-x_{d-1})}
\Big[ (y_{d-2}-1) + (y_{d-3}-1) M_1 + (y_{d-4}-1) M_1 M_2 + \ldots \Big]
\leq \frac{A'}{1+y_1-x_1}.$$
A {\em vertical section} in $Q$ is a set of the form $V_{j_1,\ldots,j_{d-2}} \cap Q$,
where \esp $V_{j_1,\ldots,j_{d-2}} := \{(j_1,\ldots,j_{d-2},i) \!: i \in \mathbb{Z} \}$.

\vspace{0.3cm}

\noindent{\bf Very good points and levels.} Assume now we are given a
positive function $\ell$ defined on $\mathbb{Z}^{d-1}$. For each
parallelepiped \esp $Q' \subset \mathbb{Z}^{d-1}$, \esp denote 
$$L_{Q'} \esp := \!
\sum_{(i_1,\ldots,i_{d-1}) \in Q'} \ell(i_1,\ldots,i_{d-1}),$$
and let \esp $\langle \ell_{Q'}\rangle$ \esp be the value above divided 
by the cardinality of $Q'$.

Given $\lambda \geq 1$, we say that a level $H_i \cap Q$ with associated sequence
(\ref{asoc-sequence}) is {\em fully $\lambda$-good} if for all $k \geq 1$,
$$\langle \ell_{Q_{m_1,\ldots,m_{k}}} \rangle \leq \lambda \langle \ell_{Q} \rangle.$$
Notice that the proportion of fully $\lambda$-good levels
is larger than \esp $(1 - \frac{d-2}{\lambda})$. Analogously,
we say that the point $p := (j_1,\ldots,j_{d-2},i_{d-1})
= V_{j_1,\ldots,j_{d-2}} \cap H_{i_{d-1}}$ is
fully $\lambda$-good whenever for all $k \geq 1$,
\begin{equation}\label{la-que-sirve}
\langle \ell_{Q_{m_1,\ldots,m_{k}} \cap V_{j_1,\ldots,j_{d-2}}} \rangle 
\leq \lambda \langle \ell_Q \rangle.
\end{equation}
For each $\lambda' \geq 1$, the proportion of fully $\lambda \lambda'$-good points
in any fully $\lambda$-good level is larger than $(1 - \frac{d-2}{\lambda} -
\frac{d-2}{\lambda'})$.

\vspace{0.2cm}

\noindent{\bf Reaching points from very good points along good vertical sections.}
A {\em segment} in $Q$ is a sequence of points for which there
exists a generator \esp $f_{i,j} \in N_d$ \esp such that each point
is obtained from the preceding one by the action of either $f_{i,j}$ or
its inverse. Such a segment $\gamma$ is said to be {\em horizontal} if
this generator is $f_{2,1}$ and $\gamma$ contains $(1+y_1-x_1)$ points.
The segment is said to be {\em vertical} if the generator involved is
one among \esp $f_{d,1},\ldots,f_{d,d-1}$ \esp (with no hypothesis
on the number of points).

Given $\lambda \geq 1$, \esp say that a segment
$\gamma$ in \esp $Q$ \esp is {\em $\lambda$-good} if
$$\langle \ell_{\gamma} \rangle \leq \lambda \langle \ell_{Q}\rangle.$$

\vspace{0.02cm}

\begin{lem} \label{vamos}
{\em Given $0 \!<\! \kappa \!<\! 1$ and $\mu \geq 1$, there exist constants 
$\lambda = \lambda_1 (\kappa, \mu, A, d)$ and $D' > 0$ such that the following 
holds: If $p := (j_1,\ldots,j_{d-2},i) = V_{j_1,\ldots,j_{d-2}} \cap H_i$ is a fully
$\mu$-good point in a $A$-round parallelepiped $Q$ such that the level
$H_i \cap Q$ is admissible, then starting from $p$ one can reach a proportion 
at least $\kappa$ of the points in $V_{j_1,\ldots,j_{d-2}} \cap Q$ by concatenating 
no more than $d-2$ vertical segments that are $\lambda$-good and have no more 
than $D' (1 + y_{1,n} - x_{1,n})$ points.}
\end{lem}

\noindent{\bf Proof.} Starting from $p$ and using $f_{d,1}^{\pm 1}$, one
can reach all points in $Q_{m_1,\ldots,m_{d-2}} \cap V_{j_1,\ldots,j_{d-2}}$ \esp
(via the segment \esp $\gamma := Q_{m_1,\ldots,m_{d-2}} \cap V_{j_1,\ldots,j_{d-2}}$).
\esp Using (\ref{la-que-sirve}) with $k = d-2$ together with (\ref{efe-de}) and
(\ref{general-est}), we obtain
\begin{eqnarray*}
\sum_{(i_1,\ldots,i_{d-1}) \in \gamma} \ell (i_1,\ldots,i_{d-1})
&=& L_{Q_{m_1,\ldots,m_{d-2}} \cap V_{j_1,\ldots,j_{d-2}}}\\
&\leq&
\frac{\mu \esp L_{Q}}{(M_{d-2} - 1) \cdots (M_1 - 1) \prod_{j=1}^{d-2} (1+y_j-x_j)}\\
&\leq&
\frac{\mu \esp 2^{d-2} A^{3d-6}}{(1+y_1-x_1)^{d-2 + \frac{(d-2)(d-1)}{2}}} L_Q
\esp\esp\esp\esp\esp
= \esp\esp\esp\esp\esp
\frac{\mu \esp 2^{d-2} A^{3d-6}}{(1+y_1-x_1)^{\frac{1}{\alpha(d)}-1}} L_Q.
\end{eqnarray*}
Hence, 
$$\langle \ell_{\gamma} \rangle \leq \mu 2^{d-2} A^{d-4} \langle \ell_Q \rangle.$$
Now, the action of $f_{d,2}$ divides $Q_{m_1,\ldots,m_{d-3}} \cap V_{j_1,\ldots,j_{d-2}}$ 
into $j_1$ segments. Given $\lambda' \geq 1$, for a proportion larger than
$(1-\frac{1}{\lambda'})$ of these segments $\gamma$, we have
\begin{eqnarray*}
\sum_{(i_1,\ldots,i_{d-1}) \in \gamma} \ell (i_1,\ldots,i_{d-1})
&\leq&
\frac{\lambda'}{j_1} L_{Q_{m_1,\ldots,m_{d-3}} \cap V_{j_1,\ldots,j_{d-2}}}\\
&\leq&
\frac{A \lambda'}{(1+y_1-x_1)} \cdot \frac{\mu L_Q}{(M_{d-3} - 1) \cdots (M_1 - 1)
\prod_{j=1}^{d-2} (1+y_j-x_j)}\\
&\leq&
\frac{\mu \esp 2^{d-3} A^{3d-7} \lambda'}{(1+y_1-x_1)^{\frac{1}{\alpha(d)} -1}} L_Q.
\end{eqnarray*}
By concatenating these segments with the previous ones, we may reach
from $p$ a proportion larger than $(1-\frac{2}{\lambda'})$ of the points of
$Q_{m_1,\ldots,m_{d-3}} \cap V_{j_1,\ldots,j_{d-2}}$.\footnote{The extra
factor $2$ comes from that the segments in consideration may
differ in number of points in $Q_{m_1,\ldots,m_{d-3}}$ by 1.}
Similarly, the action of $f_{d,3}$ divides $Q_{m_1,\ldots,m_{d-4}}$ into
$j_2$ paths; from these, a proportion larger than $(1 - \frac{2}{\lambda'})$
is $\lambda''$-good for $\lambda'' := \mu 2^{d-4} A^{3d-8} \lambda'$.
By concatenating these paths to the preceding ones, we may reach from $p$
a proportion larger than $(1-\frac{2}{\lambda'})$ of the points in
$Q_{m_1,\ldots,m_{d-4}} \cap V_{j_1,\ldots,j_{d-2}}$.

Continuing this procedure and choosing appropriately $\lambda'$, the 
concatenation property follows. Moreover, it is clear from the construction that 
the claim concerning the cardinality of each of the $\lambda$-good segments 
holds for a certain constant $D' = D'(A)$. We leave the details to the reader. $\hfill\square$

\vspace{0.3cm}

\noindent{\bf Concatenating sequences along finitely many parallelepipeds.} We 
let $\mathcal{F}_d$ be the family of finite sequences $Q^1,\ldots,Q^{d-1}$ 
of parallelepipeds in $\mathbb{Z}^{d-1}$ such that if \esp 
$Q^j = \prod_{i=1}^{d-1} [x_{i,j},y_{i,j}]$, \esp then
$$Q^{j+1}
= \cdots \times [x',y_{j,j}] \times [x_{j+1,j},y'] \times \cdots,$$
where \esp $x' > x_{j,j}$, \esp $y' > y_{j+1,j}$ \esp
(the dots mean that the corresponding entries remain untouched).
Given $A \geq 1$, we denote by $\mathcal{F}_{d,A}$ the subfamily of all sequences
made of $A$-round parallelepipeds.

Given $\mu \geq 1$, say that an horizontal segment \esp $\gamma :=
\{(i,j_2,\ldots,j_{d-1}): i \in [x_{1,1},y_{1,1}] \}$ \esp in $Q^1$ is
fully $\mu$-good with respect to $Q^1,\ldots,Q^{d-1}$ if the following
properties hold:

\vspace{0.2cm}

\noindent -- The level \esp $H_{j_{d-1}} \cap Q^{d-1}$ \esp
is admissible and fully $\mu$-good in $Q^{d-1}$.

\vspace{0.2cm}

\noindent -- The level \esp $\{ (i_1,\ldots, i_{d-3}, j_{d-2}, j_{d-1})
\!: i_k \in \mathbb{Z}\} \cap Q^{d-2}$ \esp is admissible and fully $\mu$-good
in the parallelepiped $Q^{d-2} \cap \{ (i_1,\ldots,i_{d-2}, j_{d-1}) \! :
i_k \in \mathbb{Z}\}$ (where the last intersection is understood as being
contained in \esp $\mathbb{Z}^{d-2} \sim \mathbb{Z}^{d-2} \times \{ j_{d-1} \}$).

\noindent $\vdots$

\noindent -- The level \esp
$\{ (i_1,i_2, j_3, \ldots, j_{d-2}, j_{d-1}) \!: i_k \in \mathbb{Z} \}
\cap Q^{3}$ \esp is admissible and fully $\mu$-good in the parallelepiped
$Q^{3} \cap \{ (i_1, i_2, i_3, j_4,\ldots,j_{d-1}) \! : i_k \in \mathbb{Z}\}$
(where the last intersection is understood as being contained in \esp
$\mathbb{Z}^{3} \sim \mathbb{Z}^{3} \times \{ (j_4, \ldots, j_{d-1}) \}$).

\vspace{0.2cm}

\noindent Notice that for a certain constant $A'' = A''(A,d)$, in a
proportion larger than $(1 - \frac{d-3}{\mu} - \frac{A''}{1+y_1-x_1})$,
horizontal segments in $Q^1$ are fully $\lambda$-good.

\vspace{0.15cm}

Given a sequence \esp $Q^1,\ldots,Q^{d-1}$ \esp in $\mathcal{F}_d$, a
{\em concatenating sequence from $Q^1$ to $Q^{d-1}$} is a sequence of
segments \esp $\gamma^1, \ldots, \gamma^{k}$ \esp  such that:

\begin{itemize}
\item \noindent Each $\gamma^i$ is a segment in one of the $Q^{j}$'s.

\item \noindent Each $\gamma^i$ intersects $\gamma^{i+1}$.

\item \noindent The segment $\gamma^1$ is horizontal in $Q^1$,
whereas $\gamma^{k}$ is vertical in $Q^{d-1}$.

\end{itemize}

\noindent We say that such a sequence is $\lambda$-good for $\lambda \geq 1$
if each of its segments is $\lambda$-good in one of the $Q^j$'s containing it.

\vspace{0.1cm}

\begin{lem} \label{vamos-again}
{\em Given $A > 0$, $\mu \geq 1$, and $0 \!<\! \kappa \!<\! 1$, 
there exists \esp $\lambda = \lambda_2 (\kappa,\mu,A,d)$ \esp 
such that the following holds: Let \esp $Q^1,\ldots,Q^{d-1}$ \esp be a 
sequence in $\mathcal{F}_d(A)$ and $\gamma := \{(i,j_2,\ldots,j_{d-1}) \!: i
\in [x_{1,1},y_{1,1}]\}$ a fully $\mu$-good horizontal segment for this
sequence. Then one can reach a proportion at least $\kappa$ of the points
in $Q^{d-1}$ via a $\lambda$-good concatenating sequence from $Q^1$ to $Q^{d-1}$
that starts with $\gamma^1 \!:=\! \gamma$ and is formed by no more than $K_d$
segments for a certain constant $K_d \geq 1$.}
\end{lem}

\noindent{\bf Proof.} We proceed by induction. The argument for $d = 3$ is
similar to that of the general case. It also corresponds to a more accurate
quantitative version of that given in the previous section. For this reason,
we leave it as an exercise to the reader.

Assume that the claim holds up to $d$, and let us consider the case of a sequence
\esp $Q^1,\ldots,Q^{d}$ \esp in $\mathcal{F}_{d+1}(A)$. The inductive hypothesis
applies to the sequence \esp $Q^1 \cap H_{i_{d}}, \ldots, Q^{d-1} \cap H_{i_{d}}$,
\esp where each of these intersections is understood as a parallelepiped in
$\mathbb{Z}^{d-1}$. Indeed, the definition above is made so that $\gamma$ is
also fully $\mu$-good with respect to this sequence. Accordingly, if we fix \esp
$0 \!<\! \kappa' \!<\! 1$, \esp then starting with $\gamma^{1} := \gamma$ and using
no more than $K_{d}$ segments that are $\lambda_2 (\kappa', \mu, A, d)$-good, we may
reach a proportion larger than $\kappa'$ of the points in $H_{i_{d}} \cap Q^{d-1}$.
By (\ref{efe-de}), these correspond to a proportion larger than $\kappa' +1/A^2 -1$
of the points in $H_{i_d} \cap Q^d$. This last level is fully $\mu$-good in $Q^{d}$,
hence in a proportion larger than $(1- (d-2)(1-\kappa'))$, its points are fully
$\mu/(1-\kappa')$-good. By Lemma \ref{vamos}, every such a point may reach
a proportion at least $\kappa'$ of the points in its vertical set in
$Q^{d}$ by concatenating no more than $d$ vertical segments that are
$\lambda_1(\kappa,\mu/(1-\kappa'),A,d+1)$-good. Therefore, the concatenation
of these two sequences of segments allows reaching a proportion of points
in $Q^d$ larger than
$$1-[(1-(\kappa' + 1/A^2 - 1)) + (1-\kappa')] = 2\kappa' + 1/A^2 - 2.$$
Choosing $\kappa'$ appropriately (very close to 1),
this allows concluding the proof for \esp
$K_{d+1} := K_d + d$ \esp and \esp $\lambda_2 (\kappa,\mu,A,d+1) := \max
\{\lambda_2 (\kappa', \mu, A, d), \lambda_1(\kappa,\mu/(1-\kappa'),A,d+1)\}$.
$\hfill\square$

\vspace{0.3cm}

\noindent{\bf Proof of Theorem C.} We come back to the sequence of parallelepipeds
$(Q(n))$ introduced at the beginning of \S \ref{in-general}. By the asymptotics of
the lengths of the their sides, there exists a constant $A=A_d$ such that the
following holds: For each $l \geq 0$, the finite sequence \esp $Q^{l(d-1)+1}, Q^{l(d-1)+2},
\ldots, Q^{l(d-1)+d-1}$ \esp belongs to the family $\mathcal{F}_d (A)$.

Fix $\mu \geq 1$ such that for all $n$ bigger than or equal to a certain fixed $N_0$,
$$1 - \frac{d-3}{\mu} - \frac{A''(A,d)}{1+y_{1,n(d-1)}-x_{1,n(d-1)+1}} > \frac{1}{2}.$$
Then more than a half of the horizontal segments of $Q(n(d-1)+1)$ are fully $\mu$-good.
We may hence fix $\lambda' \geq 1$ so that horizontal segments are not only
fully $\mu$-good but also $\lambda'$-good in $Q(n(d-1)+1)$ in a proportion
larger than $1/2$. By the preceding lemma, starting with any fully
$\mu$-good horizontal segment in $Q((d-1) N_0+1)$, we may find an
infinite concatenating sequence of $\lambda$-good segments for
$$\lambda := \max\{ \lambda', \lambda_2 (1/2,\mu,A_d,d) \}.$$
Moreover, according to Lemma \ref{vamos}, each of these segments contained 
in $Q(n)$ have no more than \esp $D' (1+y_{1,n} - x_{1,n})$ \esp points. 

Modulo slightly changing the constant $\lambda$ above, we may actually
assume that the sequence begins at $Q(0)$. 
We claim that the sequence that remains after cutting along concatenation points
satisfies the properties from \S \ref{in-general}. The first two ones are obvious,
while the third one easily follows from that the segments lying in a parallelepiped
$Q(n)$ must communicate between $Q(n-1)$ and $Q(n+1)$, and these two last parallelepipeds
are at distance which is comparable with the length of one of the sides of $Q(n)$. To
conclude the proof, we need to check the appropriate version of (\ref{it-works}). To do
this, just notice that if $\gamma$ is $\lambda$-good in $Q(n)$ and contains at least
$4^{\frac{n}{d-1}}/D$ points, then by H\"older's inequality,
\begin{eqnarray*}
\sum_{(i_1,\ldots,i_{d-1}) \in \gamma} \ell(i_1,\ldots,i_{d-1})^{\alpha}
&\leq& \left( \sum_{(i_1,\ldots,i_{d-1}) \in \gamma}
\ell(i_1,\ldots,i_{d-1}) \right)^{\alpha} |\gamma|^{1-\alpha}\\
&\leq& \left( \frac{\lambda L_Q}{(1+y_{1,n}-x_{1,n})^{\frac{1}{\alpha}-1}} \right)^{\alpha}
\big( D' (1+y_{1,n}-x_{1,n} \big))^{1-\alpha}\\
&=& B L_Q^{\alpha},
\end{eqnarray*}
where the last equality defines $B$. $\hfill\square$


\vspace{0.35cm}

\noindent{\bf Acknowledgments.} It is a pleasure to thank Amie Wilkinson and Benson Farb
for their invitation to Chicago-Northwestern, where the main argument of this work was
conceived.

The author was funded by the Fondecyt Research Project 1100536.


\begin{small}

\vspace{0.1cm}

\noindent Andr\'es Navas\\

\noindent Univ. de Santiago de Chile\\

\noindent Alameda 3363, Estaci\'on Central, Santiago, Chile\\

\noindent Email: andres.navas@usach.cl\\

\end{small}

\end{document}